\documentclass[final,hidelinks,onefignum,onetabnum]{siamart251216}
\usepackage{amsmath}
\usepackage{color,graphicx}
\usepackage{amssymb}
\usepackage[noend]{algpseudocode}
\usepackage{algorithm}
\usepackage{mathtools}
\usepackage{placeins}
\usepackage{pgfplots}
\usepackage{subcaption}
\usepackage[normalem]{ulem}
\usepackage{comment}
\usepackage{booktabs}
\usetikzlibrary{calc}

\newcommand{\lumo}{\lambda_\textrm{lumo}}
\newcommand{\homo}{\lambda_\textrm{homo}}
\newcommand{\lumoinner}{\lumo^{\text{in}}}
\newcommand{\lumoouter}{\lumo^{\text{out}}}
\newcommand{\lumoacc}{\lumo^{\text{acc}}}
\newcommand{\homoinner}{\homo^{\text{in}}}
\newcommand{\homoouter}{\homo^{\text{out}}}
\newcommand{\homoacc}{\homo^{\text{acc}}}

\headers{Recursive expansion of the matrix step function}{Emanuel H. Rubensson, Elias Jarlebring, and Gustaf Lorentzon}

\title{Recursive expansion of the matrix step function using polynomials of degree eight
  \thanks{\textbf{Funding:} This work was supported by the Swedish strategic research programme eSSENCE.}}
\author{Emanuel H. Rubensson\thanks{Division of Scientific Computing, Department of Information Technology, Uppsala University, Box 337, SE-751 05 Uppsala, Sweden (\email{emanuel.rubensson@it.uu.se})}
  \and
  Elias Jarlebring\thanks{Department of Mathematics, KTH Royal Institute of Technology, Lindstedtsvägen 25, SE-100 44 Stockholm, Sweden}
  \and
  Gustaf Lorentzon\footnotemark[3]}
  \ifpdf
\hypersetup{
  pdftitle={Recursive expansion of the matrix step function using polynomials of degree eight},
  pdfauthor={Emanuel H. Rubensson, Elias Jarlebring, and Gustaf Lorentzon}
}
\fi
\allowdisplaybreaks[4]
\begin{document}
\maketitle
\begin{abstract}
We consider the problem of efficiently computing the matrix step function of a large dense symmetric matrix.  To this end, we introduce a recursive polynomial expansion method in which a composite polynomial of high degree is built recursively from component polynomials of degree eight.
The component polynomial used in each iteration is designed to achieve
strong amplification of the spectral gap across the step while
favorably positioning the updated gap for subsequent iterations.  A
key ingredient is a novel evaluation scheme for arbitrary matrix
polynomials of degree exactly eight requiring only three matrix-matrix
multiplications and three matrices in memory.  This scheme makes
available a substantially larger class of component polynomials than
previously possible within a three-multiplication budget, thereby
expanding the class of composite polynomials that can be generated.
Together with our polynomial selection strategy, this leads to a
significant and consistent reduction in the number of matrix-matrix
multiplications required to compute the matrix step function compared
to existing recursive expansion methods.
\end{abstract}



\section{Introduction}
Consider a large dense symmetric matrix $X\in\mathbb{R}^{n\times n}$
whose eigenvalues lie in the union of two non-overlapping given
intervals, the left interval $I_L$ and the right interval $I_R$, i.e.,
$\lambda(X)\subset I_L\cup I_R$. We present a method for the
computation of the matrix step function for this matrix, where the
step is located in the gap.  We want to compute
\begin{align}
    D := f_\mu(X),  
\end{align}
where the parameter $\mu\in\mathbb{R}$ is the location of the step and the scalar function $f_\mu$ is defined by  
\begin{align}\label{eq:stepfuncdef}
  f_\mu(x) :=
  \begin{cases}
    0 & \textrm{if } x \leq \mu, \\
    1 & \textrm{otherwise,}
  \end{cases}
\end{align}
which we generalize to a matrix function in the standard
sense~\cite{Higham_functions_of_matrices_2008}.  The matrix $D$ is the
spectral projector onto the invariant subspace of $X$ associated with
eigenvalues greater than $\mu$.  In this paper we consider methods
that require only basic matrix operations: addition, multiplication
with a scalar, and non-scalar multiplication (matrix-matrix
multiplication), which are all well supported and parallelized on
modern high-performance computing architectures.

Our driving application stems from electronic structure calculations,
although the method and results are more generally applicable. For
consistency with the literature, we will use some of the notation from
that field.  We will use homo (highest occupied molecular orbital) and
lumo (lowest unoccupied molecular orbital) to denote the inner
boundaries of the two intervals, i.e.,
\begin{equation}
    I_L \cup I_R=[0,\lumo]\cup [\homo,1]
\end{equation}
where we additionally assumed that the problem is normalized to the
interval $[0,1]$. The normalization assumption can be made without
loss of generality in practice, if bounds on the lowest and highest
eigenvalue of $X$ are known.  The gap will be denoted
$\xi:=\lambda_{\rm homo}-\lambda_{\rm lumo}$ and referred to as the homo-lumo gap.  Since the solution to
the problem is independent of $\mu$, as long as $\mu \in (\lambda_{\rm
  lumo}, \lambda_{\rm homo})$, the problem is completely defined by
$\lambda_{\rm lumo}$ and $\lambda_{\rm homo}$; consequently our
algorithm will not involve $\mu$.

Our approach to compute $f_\mu(X)$ is based on a recursion where each
iteration $i$ involves the application of a polynomial $p_i$, called a
component polynomial, see
Algorithm~\ref{alg:general_recursive_expansion}.
\begin{algorithm}
  \caption{General recursive polynomial expansion method}
  \label{alg:general_recursive_expansion}
  \begin{algorithmic}
    \State $X_0=X$
    \For{$i=1,2,\dots$}
    \State $X_i = p_i(X_{i-1})$
    \EndFor
  \end{algorithmic}
\end{algorithm}
In this method, which we refer to as a recursive polynomial expansion,
$p_i$ are low-degree polynomials chosen so that the composite
polynomial produced by the recursive sequence converges to the target
step function:
\begin{align}\label{eq:composite_approx_problem}
  \lim_{i\rightarrow \infty} p_i(p_{i-1}(\dots p_2(p_1(x)) \dots )) = f_{\mu}(x).
\end{align}
This construction generates a polynomial of high degree with few
non-scalar multiplications. For example, the use of second-degree
component polynomials gives a composite polynomial degree $2^i$ after
$i$ iterations, or more than one million after only 20 non-scalar
multiplications.

There are various strategies for choosing the component polynomials, and
we give an overview of previous work in Section~\ref{sec:background}. One approach 
starts from a minimax perspective, where, in each iteration, one
attempts to minimize the maximum deviation from the target step
function,
\begin{align}
  \max_{x\in I_L \cup I_R} |p_i(\dots p_1(x)\dots)-f_\mu(x)|.
\end{align}
Another approach targets the asymptotic convergence toward
idempotency, where the deviation from idempotency is measured by the
idempotency error
\begin{align}
  \|X_i-X_i^2\|_2.
\end{align}
These two characterizations of convergence can be useful for the
selection of component polynomials but offer only incomplete
descriptions of the full convergence process.  In fact, during an
initial phase of the expansion, some of the most efficient methods
reduce neither the maximum deviation from the step function nor the
idempotency error.

\begin{figure}
    \centering
    \resizebox{\textwidth}{!}
    {%
        \begin{tikzpicture}[x=\textwidth, y=0.4\textwidth]
  \def\axislen{1}
  \def\vsep{0.4}
  \def\circleradius{5pt}
  \def\crosssize{8pt}

  \tikzset
    {
        icon node/.style={
            minimum width=2*\circleradius,
            minimum height=2*\circleradius,
            inner sep=0pt
        },
        text node/.style={
            minimum width=2*\circleradius,
            minimum height=2*\circleradius,
            inner sep=0pt
        }
    }
    \newcommand{\circleicon}{\tikz\draw[blue!70, line width=1.5pt] (0,0) circle (\circleradius);}
    \newcommand{\crossicon}
    {%
        \tikz{
            \draw[red!70, line width=1. 5pt] (-\crosssize/2,  \crosssize/2) -- ( \crosssize/2, -\crosssize/2);
            \draw[red!70, line width=1.5pt] ( \crosssize/2,  \crosssize/2) -- (-\crosssize/2, -\crosssize/2);
        }%
    }
    
    \newcommand{\axisval}[2]
    {
        \draw[line width=1.5pt] #1 ++(0, \circleradius) -- ++(0, -2*\circleradius);
        \node[anchor=south, yshift=+\circleradius] at #1 {#2};
    }
  \newcommand{\drawmarker}[3]
    {
        \draw[line width=1.5pt] (#1, #2) ++(0, \circleradius) -- ++(0, -2*\circleradius);
        \node[anchor=north, yshift=-\circleradius] at (#1, #2) {#3};
    }
    \newcommand{\drawmarkerup}[3]
    {
        \draw[line width=1.5pt] (#1, #2) ++(0, \circleradius) -- ++(0, -2*\circleradius);
        \node[anchor=south, yshift=+\circleradius] at (#1, #2) {#3};
    }
    \newcommand{\drawmarkergray}[3]
    {
        \draw[gray!70, line width=1.5pt] (#1, #2) ++(0, \circleradius) -- ++(0, -2*\circleradius);
        \node[anchor=north, yshift=-\circleradius] at (#1, #2) {#3};
    }
    \newcommand{\drawcross}[2]{
      \draw[red!70, line width=1.5pt] (#1, #2) ++(-\crosssize/2,  \crosssize/2) -- ++( \crosssize, -\crosssize);
      \draw[red!70, line width=1.5pt] (#1, #2) ++( \crosssize/2,  \crosssize/2) -- ++(-\crosssize, -\crosssize);
    }
    \newcommand{\drawcircle}[2]{
      \draw[blue!70, line width=1.5pt] (#1, #2) circle (\circleradius);
    }

    \pgfplotstableread[col sep=comma]{occ-data.csv}\occdata
    \pgfplotstableread[col sep=comma]{unocc-data.csv}\unoccdata
    \pgfplotstableread[col sep=comma]{gap-data.csv}\gapdata

    \foreach \k in {0,1,2,3}
    {
        \pgfmathsetmacro{\y}{-\k * \vsep}
        \coordinate (axisstart-\k)  at (-0.05, \y);
        \coordinate (axisend-\k)    at (\axislen+0.05, \y);
        \coordinate (ax1-\k-0)      at (0, \y);
        \coordinate (ax1-\k-1)      at (1, \y);
        \coordinate (ax1-\k-half)   at (0.5, \y);

        \draw[->, line width=1.5pt] (axisstart-\k) -- (axisend-\k);
        \axisval{(ax1-\k-0)}{$0$}
        \axisval{(ax1-\k-1)}{$1$}
        \axisval{(ax1-\k-half)}{$0.5$}
    }

    \foreach \k in {0,1,2,3}
    {
        \pgfmathsetmacro{\y}{-\k * \vsep}
        
        \pgfplotstableforeachcolumnelement{k\k}\of\gapdata\as\val{
            \drawmarkergray{\val}{\y}{$\mu^{(\k)}$}
        }   
        
        %
        \def\idx{0}        \pgfplotstableforeachcolumnelement{k\k}\of\occdata\as\val{
            \pgfmathtruncatemacro{\idx}{\idx+1}
            \node[icon node] (occ-{\k}-{\idx}) at (\val, \y) {\circleicon};
        }
        %
        \def\idx{0}
        \pgfplotstableforeachcolumnelement{k\k}\of\unoccdata\as\val{
            \pgfmathtruncatemacro{\idx}{\idx+1}
            \drawcross{\val}{\y}
            \node[icon node] (unocc-{\k}-{\idx}) at (\val, \y) {\crossicon};
        }
        %

        \node[anchor=north, yshift = -3pt] (lumosymb) at (occ-{\k}-{4}) {$\lumo^{(\k)}$};
        \node (gapleft) at ($(lumosymb.south) + (0, -6pt)$) {};
        \draw[gray!70, line width=1.5pt] (gapleft) -- ++(0, 5pt) -- ++(0, -10pt);

        \node[anchor=north, yshift = -3pt] (homosymb) at (unocc-{\k}-{1}) {$\homo^{(\k)}$};
        \node (gapright) at ($(homosymb.south) + (0, -6pt)$) {};
        \draw[gray!70, line width=1.5pt] (gapright) -- ++(0, 5pt) -- ++(0, -10pt);

        \node (gapcenter) at ($(gapleft)!0.5!(gapright)$) {$\xi^{(\k)}$};

        \draw[<-] (gapleft) -- (gapcenter.west) {};
        \draw[->] (gapcenter.east) -- (gapright) {}; 
    }
    \foreach \k in {0,1,2,3}
    {
        \ifnum\k=0
            \node[anchor = east] at (axisstart-\k) {$X_0$};
        \fi
        \pgfmathtruncatemacro{\kminusone}{\k - 1}
        \ifnum\k>0
            \node[left] at (axisstart-\k) {$X_\k = p_\k(X_{\kminusone})$};
        \fi
    }

    
    \node[icon node] (legendcircle) at (0, -4.2*\vsep) {\circleicon};
    \node[anchor= west] (legend1) at (legendcircle.east) {Left Partition of Eigenvalues of $X_{i}$};
    
    \node[icon node, anchor=west] (legendcross) at ({0.5,0} |- legend1) {\crossicon};
    \node[anchor=west] (legend2) at (legendcross.east) {Right Partition of Eigenvalues of $X_i$};
\end{tikzpicture}
    }
    \caption{ Eigenvalue spectrum vs iteration in a polynomial recursion. The polynomials $p_1(x) = x^2$, $p_2(x) = 2x-x^2$ and $p_3(x) = 2x-x^2$ were used. These polynomials were chosen in accordance with the SP2 scheme described in Section~\ref{sec:recursive-methods}.}
    \label{fig:eigenevolution}
\end{figure}
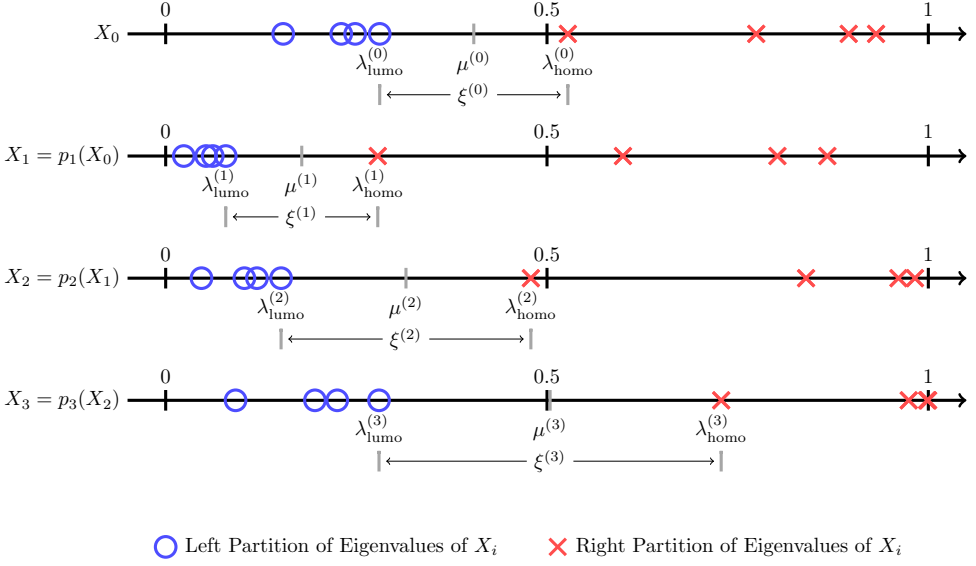

In the initial phase of the expansion it is more appropriate to consider
how fast the expansion amplifies the gap.
In each step of the recursion the left and right intervals are updated
by the transformation
\begin{equation}
I_L^{(i)}=\left[0,\;\max_{x\in I_L^{(i-1)}}(p_{i}(x))\right],\;\; 
I_R^{(i)}=\left[\min_{x\in I_R^{(i-1)}}(p_{i}(x)),\;1\right],\;\; i = 1,2,\dots,
\end{equation}
under a normalization condition of the component polynomials which
does not restrict generality.  Each step leads to a new problem of the
same type but with a new location and size of the gap, see Fig.~\ref{fig:eigenevolution}.  We may define
the condition number of the problem in iteration $i$
as in~\cite{Rubensson_2008} where the following identity between the condition number and the inverse of the gap is shown:
\begin{align}\label{eq:condition_number}
  \lim_{h\rightarrow 0} \sup_{E:\|E\|_2 =
    1}\frac{\|f_{\mu^{(i)}}(X_i+hE)-f_{\mu^{(i)}}(X_i)\|_2}{h} = \frac{1}{\xi^{(i)}}.
\end{align}
This means that a recursive expansion is forward stable if $\xi^{(i)}\geq \xi^{(0)},\ i=1,2,\dots$. 
In the first phase of the expansion the condition number is reduced
and we therefore call it \emph{the conditioning phase}.  When the
condition number comes close to 1, the idempotency error becomes
relevant and the asymptotic order of convergence becomes important. We
refer to this final phase of the expansion as \emph{the purification
phase}. We adopt this terminology from~\cite{Rubensson_2014}.

In this work we aim to select component polynomials that, for a given
number of non-scalar multiplications, reduce the condition number as
quickly as possible in the conditioning phase and achieve as high an
order of convergence as possible in the purification phase.
The novelties are
\begin{itemize}
\item a recursive polynomial expansion method based on a new family of
  degree-eight component polynomials, where the polynomials are
  constructed to satisfy equioscillation conditions on each of the
  left and right intervals,
\item a scheme for the evaluation of any polynomial of
  degree exactly eight, requiring only three non-scalar multiplications and three
  matrices in memory.
\end{itemize}

Our choice of degree eight is motivated by a new approach for the
evaluation of polynomials by Sastre~\cite{Sastre_2018, Alonso_2026}, who showed that one can
compute most polynomials of degree eight with only three non-scalar
multiplications. We further develop this method in
Section~\ref{sec:eval}, to avoid complex arithmetic and algorithm
breakdown occurring for certain values of polynomial coefficients. We also provide a scheme requiring only three matrices in memory.
This means that for $d\leq 3$ we now have that $d$ multiplications are
sufficient to evaluate arbitrary polynomials strictly of degree $2^d$.
It was recently shown that corresponding schemes do not exist for
$d>3$~\cite{Jarlebring_2026}, and degree eight is therefore the highest degree for
which this full representability property holds.

While our degree-eight component polynomials are not derived
from a minimax approximation problem, their construction is
inspired by the equioscillation principle underlying
minimax polynomial approximations. For our component 
polynomials, we require a total of seven equioscillations 
over the left and right intervals. 
For any
given pair $(\lumo,\, \homo)$, the equioscillation conditions yield a
family of eight candidate polynomials characterized by the number of
oscillations that fall in the left and right intervals.  In each
iteration we select the polynomial in the family that maximizes the
gap.  A key design decision of our method is to not impose equal
oscillation amplitudes on the left and right intervals.
This freedom allows $\mu^{(i)}$ to be located arbitrarily in $[0,\,1]$
and our polynomials are therefore distinct from polynomials resulting
from a greedy minimax approach that effectively locks $\mu^{(i)}$ at
0.5.  We will show that this results in faster convergence as greater
gap amplification can be achieved when $\mu^{(i)}$ is located away
from 0.5.

\section{Background and related work}\label{sec:background}
A key computational task in electronic structure theory, based on for
example Hartree-Fock or density functional theory, is the computation
of the single-particle density matrix $D$ for a given Fock or
Kohn-Sham matrix $H$. The density matrix is given by
\begin{align}\label{eq:density_matrix}
  D = f_\mu(X_0),
\end{align}
where 
\begin{align}\label{eq:initial_transformation}
X_0 := (l_{\max}I-H)/(l_{\max}-l_{\min}) 
\end{align}
and $l_{\min}$ and $l_{\max}$ are lower and upper bounds of the
eigenspectrum of $H$, respectively.  
The density matrix is the spectral projector associated with 
eigenvalues of $X_0$ larger than $\mu$, and its trace equals 
the number of occupied states, $n_{\text{occ}} = \text{Tr}(D)$.
Note that \eqref{eq:initial_transformation} implies
that $X_0$ will have its eigenvalues on the interval $[0,1]$, and
\eqref{eq:density_matrix} therefore matches the problem of this paper.  

The matrix step function problem occurs in electronic structure
calculations with several variations. The evaluation of the matrix
step function commonly appears as one step in an outer iteration,
such as a
self-consistent field iteration or molecular dynamics simulations
requiring repeated electronic structure calculations. The step
location $\mu$ is often determined implicitly so as to yield a
prescribed number of electrons, i.e.\ $\text{Tr}(D) = n_{\text{occ}}$.
Practical implementations use both dense and sparse matrix
representations, with screening of small matrix elements
usually employed in the method class \emph{linear scaling methods}~\cite{BowlerMiyazaki_2012, Rubensson_chapter_2011}. A multitude of electronic
structure models exist, for example through different energy
functionals in density functional theory, and together with the choice
of basis set the selected model affect the spectral properties and
sparsity structure of $H$. We also note that insulating systems are
computationally more favorable than metallic systems because the
spectral gap improves the conditioning of the step function (see
\eqref{eq:condition_number}) and leads to faster decay of matrix
elements, enabling sparser matrix representations~\cite{Benzi_2013}.

The recursive expansions considered in this work stand in contrast to
serial polynomial expansions of Chebyshev type that are also rather
common for the problem above~\cite{Goedecker_1994, Goedecker_1995,
  Liang_2003}.  Serial expansions are more flexible in the sense that
any polynomial of the chosen degree can be represented, but are often
less efficient than the recursive expansions in the present case, as 
the polynomial degree that can be reached for a given number of
non-scalar multiplications is much 
lower~\cite{Rubensson_chapter_2011, Niklasson_chapter_2011}.

\subsection{Methods based on recursive expansions}\label{sec:recursive-methods}
\begin{figure}
  \centering  
  \resizebox{\textwidth}{!}
  {
      \begin{subfigure}{0.276\textwidth}
        \centering
        \includegraphics[width=\linewidth, trim=29pt 10pt 65pt 15pt,clip]{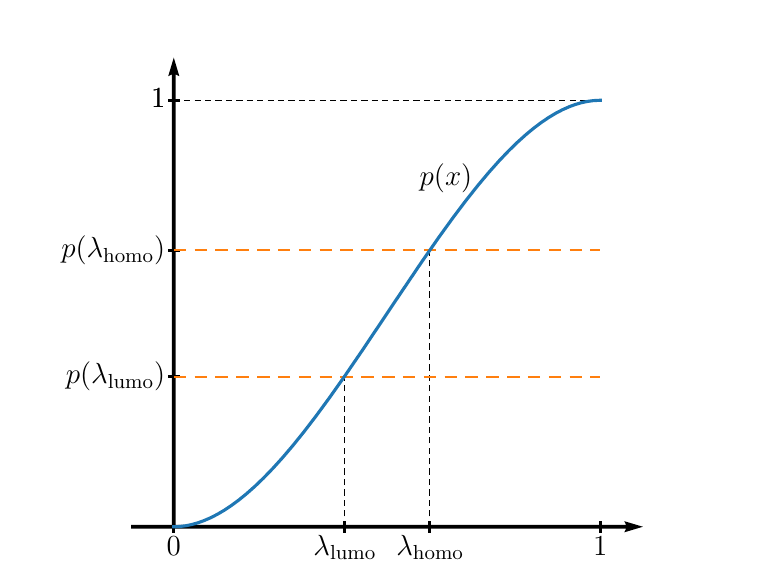}
    \caption{McWeeny \label{fig:mcweeny}}%
      \end{subfigure}
    \begin{subfigure}{0.32\textwidth}
        \centering
        \includegraphics[width=\linewidth, trim=25pt 10pt 25pt 15pt,clip]{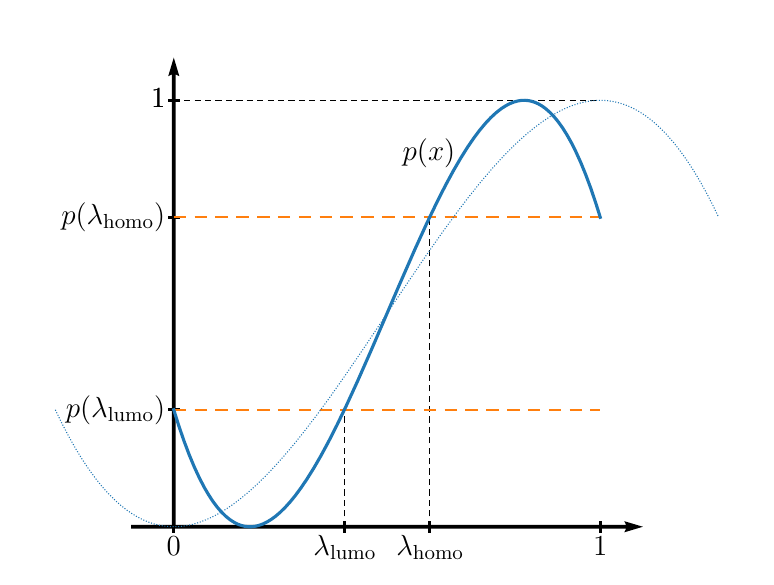}
    \caption{Accelerated McWeeny\label{fig:sp3_acc}}%
      \end{subfigure}
    \begin{subfigure}{0.276\textwidth}
        \centering
        \includegraphics[width=\linewidth, trim=29pt 10pt 68pt 15pt,clip]{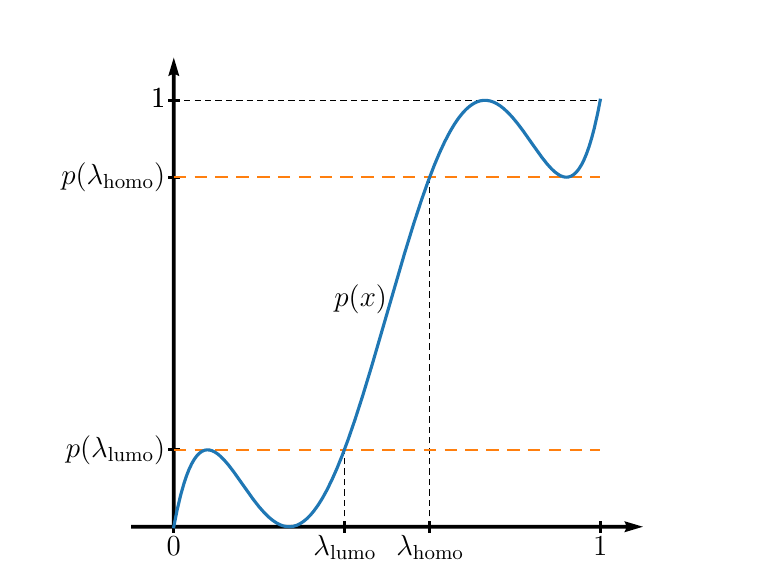}
    \caption{Su5 \label{fig:sp5_acc}}
      \end{subfigure}%
  }
  \caption{Illustration of the development from the McWeeny polynomial to accelerated recursive polynomial expansion schemes. Panel~(a) shows the McWeeny polynomial, panel~(b) its scale-and-fold accelerated variant (Su3), and panel~(c) Suryanarayana's degree-five generalization (Su5). The dotted line in panel~(b) shows the McWeeny polynomial evaluated on the stretched interval, illustrating the scale-and-fold construction.\label{fig:rec_exp}}
\end{figure}

The problem we are considering can be viewed as a variation of the
computation of the matrix sign function. The matrix sign function has
been studied for decades, including work by Denman and Beavers
\cite{denman1976matrix} and earlier work by Schulz
\cite{schulz1933iterative}. See a summary and review in Higham's monograph~\cite{Higham_functions_of_matrices_2008}.  These works
include a number of methods that fit in the recursive or iterative
structure of Algorithm~\ref{alg:general_recursive_expansion}. Several of these methods, such as for
example the well known Newton iteration, uses inverses, which is beyond our scope.
A popular method which does not require inverses is the Newton-Schulz iteration~\cite{schulz1933iterative}, based on the polynomial
$\frac{1}{2}x(3-x^2)$. After a linear transformation of this iteration
to the $[0,\,1 ]$ interval, we obtain the polynomial $3x^2-2x^3$, which
in the electronic structure literature is known as the McWeeny
polynomial~\cite{McWeeny_1956}, see Fig.~\ref{fig:mcweeny}.

Iteration with the McWeeny polynomial has been the starting point for the
development of a large number of methods targeting the matrix step
function and efficient construction of the density matrix.  An
important distinction between different methods is whether or not they
assume information about the location of $\mu$, or the $\lumo$ and
$\homo$ eigenvalues.
In a seminal paper, Palser and Manolopoulos proposed a modification to
the McWeeny iteration that automatically adjusts the step location to
give the correct number of occupied orbitals~\cite{Palser_1998}.  In
this scheme $\mu$ is allowed to move between the iterations, which
stands in contrast to schemes that use component polynomials with a
fixed point $p_i(\mu)=\mu,\ i=1,2,\dots$ (such as the McWeeny
polynomial).
The second-order spectral projection scheme (SP2) 
is a variable-$\mu$ iteration that uses a combination of the
component polynomials
$x^2$ and $2x-x^2$~\cite{Niklasson_2002}. In the original SP2 scheme, the
choice between the two 
polynomials is made using trace correction, i.e., based on whether 
$\text{Tr}[X_{i-1}]$ is above or below the desired occupation $n_{\text{occ}}$, 
thereby automatically adjusting the step location. This method outperforms
the Palser-Manolopoulos scheme when the step is away from the spectral center.
The McWeeny, Palser-Manolopoulos, and SP2 schemes all belong to a method class (see also,
e.g.,~\cite{Holas_2001, Mazziotti_2003, Niklasson_trs4_2003} ) that
uses component polynomials which are monotone on $[0,\, 1]$, have
fixed points at 0 and 1 and one or more vanishing derivatives at 0
and/or 1.

\subsection{Accelerated recursive polynomial expansions}
Faster convergence can be achieved by relaxing the monotonicity
constraints, an idea underlying a number of works exploring
accelerated recursive expansions~\cite{Rubensson_nonmonotonic_2011,
  KimJung_2011, Suryanarayana_2013, Finkelstein_2021}.  A
scale-and-fold acceleration technique was applied to the SP2 and
McWeeny polynomials to speed up
convergence~\cite{Rubensson_nonmonotonic_2011}.
The scaling in these schemes relies on reliable bounds on the $\homo$
and $\lumo$ eigenvalues to ensure both efficiency and stability.
It was shown in~\cite{Rubensson_2014} how such bounds can be obtained
as a cheap byproduct of the SP2 expansion. 
Therefore, when the step function computation is a part of an outer loop, the eigenvalue bounds can be obtained from the previous outer iteration and adjusted to account for changes in the matrix.
The scale-and-fold technique applied to the McWeeny polynomial is
illustrated in Fig.~\ref{fig:sp3_acc}. When applied recursively, this
results in an accelerated scheme.  The accelerated McWeeny scheme was
generalized to higher order by
Suryanarayana~\cite{Suryanarayana_2013}.  Suryanarayana's polynomials
have odd degree, are anti-symmetric around 0.5, such that $p(x)+p(1-x)
= 1$, and are constructed to satisfy equioscillation conditions on
each of the intervals $[0,\,\lumo]$ and $[\homo,\, 1]$, with the
homo-lumo gap centered around 0.5.
These polynomials are illustrated in Fig.~\ref{fig:sp3_acc} and
Fig.~\ref{fig:sp5_acc} for degree three and five, respectively. In the
case of degree three the Suryanarayana polynomial coincides with the
accelerated McWeeny polynomial.
We will refer to the Suryanarayana polynomial of degree X as SuX. 

In~\cite{Finkelstein_2021} the SP2 scheme was mapped to the structure of a deep
neural network and parameters of each layer were optimized leading to
a significant acceleration.
Analogous acceleration schemes have been developed for the matrix sign
function, motivated by applications in homomorphic
encryption~\cite{Cheon_2020,Lee_2022} and deep learning
optimization~\cite{jordan2024muon,Amsel_2025,Grishina_2025}.
A scheme closely related to the accelerated McWeeny iteration,
employing a heuristic scaling parameter, was used as part of a
divide-and-conquer eigensolver~\cite{Huss-Lederman_1997}.
Several of these schemes are essentially equivalent up to a linear 
transformation: the deep neural
network optimization formulation of Finkelstein et
al.~\cite{Finkelstein_2021} was shown to coincide with the
scale-and-fold accelerated SP2 expansion of
Rubensson~\cite{Rubensson_nonmonotonic_2011}, while Suryanarayana's
schemes~\cite{Suryanarayana_2013} are closely related to the
sign-function schemes proposed by Lee et al.~\cite{Lee_2022} and Amsel
et al.~\cite{Amsel_2025}.

Recursive expansions have been implemented on a variety of
computational platforms using different matrix representations. Many
of the original works focus on sparse matrix formulations aiming at
linear scaling with system size, see e.g.~\cite{Palser_1998,
  Niklasson_2002, Rubensson_nonmonotonic_2011}.  
In such sparse matrix implementations small matrix elements are
neglected, and the forward error can be controlled provided that the
norm of the error matrix introduced by those omissions can be kept
below a prescribed bound~\cite{Rubensson_2008}.
Several sparse matrix
implementations have been developed for distributed-memory
clusters~\cite{Dawson_2018, Kruchinina_2019, VandeVondele_2012}.
Dense matrix implementations have also been reported for
distributed-memory clusters~\cite{Chow_2015, Khadatkar_2023} as well
as for hardware accelerators such as GPUs~\cite{Cawkwell_2012}
(including NVIDIA tensor cores~\cite{Finkelstein_2021}) and Google
tensor processing units~\cite{Pederson_2023}. Recursive expansions
have further been used to drive graph-based electronic structure
theory~\cite{Niklasson_graph_2016, Kulichenko_graph_2025,
  Katbashev_graph_2025}.

\subsection{Efficient polynomial representation}\label{sec:efficientpolynomialrepresentation}
All recursive expansions discussed above impose structural restrictions on the polynomials that can be represented, stemming from the composite nature of the iteration and the fixed low degree of the component polynomials.  Additional constraints are often introduced to facilitate the determination of polynomial coefficients or to enforce a desired asymptotic convergence order, such as fixed points or vanishing derivatives at 0 and/or 1, or anti-symmetry around 0.5. Furthermore, certain properties may be preferred for structural or computational reasons, for example odd polynomials in sign-function iterations for the polar factor~\cite{Amsel_2025} or low component polynomial degree to limit the number of temporary matrices.

Our degree-eight component polynomials relaxes several structural restrictions of previous approaches.
Achieving degree eight with three multiplications can be compared with Suryanarayana's schemes: Su5 achieves degree five with three multiplications, and Su7 achieves degree seven with four.  Three iterations of the SP2 expansion (accelerated or not) produce a polynomial of degree eight in three multiplications, but the composite structure $p_{i+2}(p_{i+1}(p_i(x)))$ restricts the class of attainable polynomials. 
Our component polynomials therefore provide access to a
larger class of composite polynomials for a given number of 
non-scalar multiplications.

\section{Evaluation of degree-eight matrix polynomials}\label{sec:eval}
Let $b_0,\ldots,b_8\in\mathbb{R}$ be the monomial coefficients of a matrix polynomial of degree eight, i.e.,
\begin{align}\label{eq:mono8}
  p(X) = \sum_{i=0}^8 b_iX^i.
\end{align}
Sastre showed that such a degree-eight matrix polynomial can
often be evaluated using only three non-scalar multiplications~\cite{Sastre_2018}.
Sastre makes the ansatz
\begin{subequations}\label{eq:sastreeval}
\begin{align}
  y_{02}(X) =& X^2(c_4+c_3X) \\
  y_{12}(X)     =& (y_{02}(X) + d_2X^2 + d_1X)(y_{02}(X)+e_2X^2) \\
  & +e_0y_{02}(X) + f_2X^2+f_1X+f_0I\nonumber
\end{align}
\end{subequations}
and finds the nine unknown coefficients by taking $y_{12}(X) =
\sum_{i=0}^8 b_iX^i$ and equating the coefficients for each of
the matrix monomials $X^i,\, i=0,1,\dots,8$. This gives a system of
nine equations.
Unfortunately, this system has no solution for
certain values of polynomial coefficients, such as for example if
$b_8=0$ and $b_7\neq 0$ or if $b_7=b_5=0$ and $b_3\neq 0$.
Furthermore, some of the coefficients in \eqref{eq:sastreeval} are given by the solution of
quadratic equations, that have complex solutions for some polynomials. 
This leads
to the need for complex matrix-matrix multiplication even if $X$ is
real.
One such example is the polynomial $21x^8-48x^7+28x^6$, which is used later in this article.

Here, we use a similar ansatz but add one degree of freedom that is
used to avoid complex coefficients.  Our scheme requires only
$b_8\neq 0$, and always gives real coefficients.

Let
\begin{subequations}\label{eq:Smatrices}
\begin{align}
  S_0 = & I \\
  S_1 = & X \\
  S_2 = & S_1^2 \label{eq:Smatrices:mul1} \\ 
  S_3 = & S_2 (c_1 S_1 + S_2) \label{eq:Smatrices:mul2} \\
  S_4 = & (d_0 S_0 + d_1 S_1 + d_2 S_2 + S_3) (e_1 S_1 + e_2S_2 + S_3) \label{eq:Smatrices:mul3} \\
  S_5 = & f_0 S_0 + f_1 S_1 + f_2 S_2 + f_4 S_4 \label{eq:Smatrices:output}
\end{align}
\end{subequations}

By expanding $S_5(X)$ in monomials, taking $S_5(X) =
\sum_{i=0}^8 b_iX^i$, and equating the coefficients for each of the
monomials we get the equations
\begin{subequations}
    \begin{align}
      b_8 & = f_4 \label{eq:coeffA8}\\
      b_7 & = 2f_4c_1  \\
      b_6 & = f_4(c_1^2 + d_2 + e_2)  \\
      b_5 & = f_4(c_1(d_2 + e_2) + d_1 + e_1)  \label{eq:coeffA5}\\
      b_4 & = f_4(c_1(d_1 + e_1) + d_0 +d_2e_2)  \label{eq:coeffA4}\\
      b_3 & = f_4(c_1d_0 + d_1e_2 + d_2e_1)   \label{eq:coeffA3}\\
      b_2 & = (f_2 + f_4(d_0e_2 + d_1e_1))  \label{eq:coeffA2}\\
      b_1 & = (f_1 + f_4d_0e_1) \\
      b_0 & = f_0.  \label{eq:coeffI}
    \end{align}
\end{subequations}
Note that at this point we have nine equations and ten unknowns. We have designed a system which is intentionally underdetermined, and the additional unknown now allows us to build an evaluation with the desired properties. From
equations \eqref{eq:coeffA8}-\eqref{eq:coeffA5} we obtain
\begin{align}
  f_4 & = b_8,  \label{eq:coeff_deg8term}\\
  c_1 & = b_7/(2f_4), \\
  t_2 & := d_2 + e_2 = b_6/f_4 - c_1^2, \\
  t_1 & := d_1 + e_1 = b_5/f_4 - c_1t_2, \label{eq:coeff_t1}
\end{align}
where we introduced two new parameters $t_1,t_2$ for convenience.
From \eqref{eq:coeffA4} we may formulate a quadratic equation in $e_2$:
\begin{align}
  e_2^2 - t_2 e_2 + b_4/f_4 - d_0 - c_1 t_1 & = 0
\end{align}
with the roots
\begin{align}
  e_2 = \frac{1}{2}\left(t_2 \pm \sqrt{t_2^2 - 4(b_4/f_4-d_0-c_1t_1)}\right).
\end{align}
Furthermore we have from \eqref{eq:coeffA3} that
\begin{align}
  (e_2-d_2)e_1 = c_1d_0 + t_1e_2 - b_3/f_4
\end{align}
where we used that $d_1 = t_1-e_1$. We now make use of our extra
degree of freedom and add a tenth equation
\begin{align}
  t_2^2 - 4(b_4/f_4-d_0-c_1t_1) = 1.
\end{align}
Solving for $d_0$ and choosing the positive branch for $e_2$ we obtain
\begin{align}
  d_0 & = \frac{1}{4} ( 1 - t_2^2 + 4b_4/f_4 - 4c_1t_1), \label{eq:coeffd0}\\
  e_2 & = \frac{1}{2} ( t_2 + 1 ), \\
  d_2 & = t_2 - e_2 =  \frac{1}{2} ( t_2 - 1 ), \\
  e_1 & = c_1d_0 + t_1e_2 - b_3/f_4,
\end{align}
where in the last equation we used that $e_2-d_2=1$. We have thus
avoided both complex coefficients and division by zero, provided only
that $b_8\neq 0$. Finally we have from
\eqref{eq:coeffA2}-\eqref{eq:coeffI} that
\begin{align}
  f_2 & = b_2 - f_4(d_0e_2 + d_1e_1), \\
  f_1 & = b_1 - f_4d_0e_1, \\
  f_0 & = b_0. \label{eq:coeff_deg0term}
\end{align}
This shows that there exist coefficients such that \eqref{eq:Smatrices} computes $p(X)$ for a given polynomial $p$. We summarize this result in the following theorem.
\begin{theorem}
    Let $p$ be given as in \eqref{eq:mono8} with $b_8 \neq 0$. Then, $p(X)$ can be computed in 3 matrix-matrix multiplications, involving only real evaluation coefficients. 
\end{theorem}

\subsection{Memory-efficient evaluation} \label{sec:poly8_memory_efficient}
In modern high-performance computing environments, economical use of memory resources is of great importance.  In the following, we show how our algorithm for evaluation of degree-eight matrix polynomials can be employed with only three matrices stored in memory. If we choose the second multiplication to be a square of two matrices, we can double the degree without overwriting any previous information. Luckily, this can be done without loss of generality. 
However, setting up the third multiplication in the limited workspace
requires some additional work in rearranging the workspace.
The setup with the fewest number of linear combinations we managed to find is given as follows. 

Denote the workspace as $M_1$, $M_2$, $M_3$, and assume the initial state to be
\begin{align}
    M_1 \leftarrow X, \quad M_2 \text{ empty}, \quad M_3 \text{ empty}.
\end{align}
In the following, we show how to evaluate $M_1 \leftarrow p(X)$ within this workspace. The first two multiplications are computed as follows
\begin{alignat}{2}
    M_2 &\leftarrow M_1 \cdot M_1 &&= X^2,\\
  M_2 &\leftarrow M_2 + \frac{1}{2}c_1 M_1 &&= \frac{1}{2}c_1 X + X^2,\\
  M_3 &\leftarrow M_2\cdot M_2 &&= \frac{1}{4} c_1^2 X^2 + c_1 X^3 + X^4.
\end{alignat}
The workspace is then given by:
\begin{align}\label{eq:workspace_post_mul2}
    M_1 = X, \quad
    M_2 = \frac{1}{2}c_1 X + X^2, \quad
    M_3 = \frac{1}{4}c_1^2X^2 + c_1X^3 + X^4.
\end{align}
Next, we rearrange our workspace to be
\begin{subequations}\label{eq:workspace_pre_mul3}
\begin{align}
    M_1 &= f_0 I + f_1 X + f_2 X^2, \\
    M_2 &= e_1 X + e_2X^2 + c_1 X^3 + X^4, \\
    M_3 &= d_0 I + d_1 X + d_2 X^2 + c_1 X^3 + X^4,
\end{align}
\end{subequations}
using a sequence of linear combinations. Observe that $M_2$ and $M_3$ are the multiplication factors used to construct $S_4$ in \eqref{eq:Smatrices:mul3}. With this setup, we can construct the final output to be the same as in \eqref{eq:Smatrices:output}:
\begin{align}
    M_1 \leftarrow M_1 + f_4 M_2 M_3 = S_5.
\end{align}

It remains to show how we set up \eqref{eq:workspace_pre_mul3} from \eqref{eq:workspace_post_mul2}. We define some temporary scalars for convenience:
\begin{subequations}\label{eq:rscalars_def}
    \begin{align}
          r_1 &\coloneqq d_1 - \frac{1}{2}c_1(d_2 - \frac{1}{4}c_1^2), \\
          r_2 &\coloneqq d_2 - \frac{1}{4}c_1^2,    \\
          r_3 &\coloneqq e_1 - d_1 - \frac{1}{2} c_1,  \\  
          r_4 &\coloneqq f_1 - f_2(e_1-d_1).
    \end{align}
\end{subequations}
The workspace can be rearranged through the following steps
\begin{alignat}{2}
    M_3 &\leftarrow M_3 + r_1 M_1 + r_2 M_2 &&=d_1 X + d_2 X^2 + c_1 X^3 + X^4, \\
    M_2 &\leftarrow M_2 + r_3M_1 &&= (e_1 - d_1)X + X^2,\\
    M_1 &\leftarrow r_4 M_1 + f_2 M_2 + f_0 I &&= f_0 I + f_1 X + f_2 X^2,\\
    M_2 &\leftarrow M_2 + M_3 &&= e_1 X + \underbrace{(1+d_2)}_{=e_2}X^2 + c_1 X^3 + X^4,\\
    M_3 &\leftarrow M_3 + d_0 I &&= d_0 I + d_1 X + d_2 X^2 + c_1 X^3 + X^4,
\end{alignat}
after which the workspace corresponds to \eqref{eq:workspace_pre_mul3}. The process is summarized in Algorithm~\ref{alg:p8mem3}
\begin{algorithm}
\caption{Memory efficient degree-eight matrix polynomial evaluation}
\label{alg:p8mem3}
\begin{algorithmic}[1] 
\Statex \textbf{Input:} monomial coefficients $b_0,\ldots,b_8$ with $b_8 \neq0$;
memory slots $M_1,M_2,M_3,$  where $M_1=X$
\State Compute $f_0,f_1,f_2,f_4,e_1,d_0,d_1,d_2,c_1$ according to \eqref{eq:coeff_deg8term}-\eqref{eq:coeff_t1},\eqref{eq:coeffd0}-\eqref{eq:coeff_deg0term}
\State Compute $r_1,r_2,r_3,r_4$ according to \eqref{eq:rscalars_def}
\State $M_2 \leftarrow M_1M_1$
\State $M_2 \leftarrow M_2 + \frac{1}{2}c_1M_1$
\State $M_3 \leftarrow M_2 M_2$
\State $M_3 \leftarrow M_3 + r_1M_1$
\State $M_3 \leftarrow M_3 + r_2M_2$
\State $M_2 \leftarrow M_2 + r_3 M_1$   
\State $M_1 \leftarrow r_4 M_1 + f_2 M_2 + f_0 I$  
\State $M_2 \leftarrow M_2 + M_3$  
\State $M_3 \leftarrow M_3 + d_0 I$
\State $M_1 \leftarrow f_4 M_2 M_3 + M_1$
\Statex \textbf{Output:} $M_1 \leftarrow p(M_1)$
\end{algorithmic}
\end{algorithm}

\section{A family of degree-eight polynomials}\label{sec:parameterization}
We now turn to the construction of the degree-eight polynomials used in the recursive expansion. Our overall goal is to compute $f(X_0)$ using as few non-scalar multiplications as possible. In principle, the polynomials should be chosen to optimize the performance of the entire expansion, but such a global optimization is not tractable. 
Instead, we aim to construct and select in each iteration a polynomial that amplifies the homo-lumo gap as much as possible, while accounting for the location of the updated gap, as this affects the amplification achievable in subsequent steps. 
In this section, we derive a family of polynomials, referred to as the SP8 family, from which we choose one in each recursive expansion iteration.

Our polynomials are designed to satisfy equioscillatory conditions on each of the intervals $[0,\,\lumo]$ and $[\homo,\, 1]$.  For any given pair of eigenvalues $\lumo$ and $\homo$, we define a family of eight polynomials, one for each choice of $L=0,\dots,7$ stationary points in the left interval and $R=7-L$ in the right. 

To represent our polynomials we use the parameterization 
\begin{align}\label{eq:polynomial_parameterization}
    p(x) & = s_8 \sum_{k=0}^{7} (-1)^k \frac{e_k(s)}{8-k}\, x^{8-k} + s_9,
\end{align}
where 
\begin{align}
e_k(s) &= \sum_{1 \le i_1 < \cdots < i_k \le 7}
    s_{i_1} \cdots s_{i_k}, \qquad k = 0,\dots,7,
\end{align}
are the elementary symmetric polynomials in the stationary points $s_1,\dots,s_7$ with $e_0(s)=1$. Differentiating yields
\begin{align}
  p'(x) & = s_8\prod_{k=1}^{7} (x-s_k),
\end{align}
confirming that the polynomial has stationary points at $s_1,\dots,s_7$.

\begin{figure}
  \centering  
  \begin{subfigure}{0.48\textwidth}
    \centering
    \includegraphics[width=\linewidth,
      trim=40pt 10pt 60pt 15pt,clip]{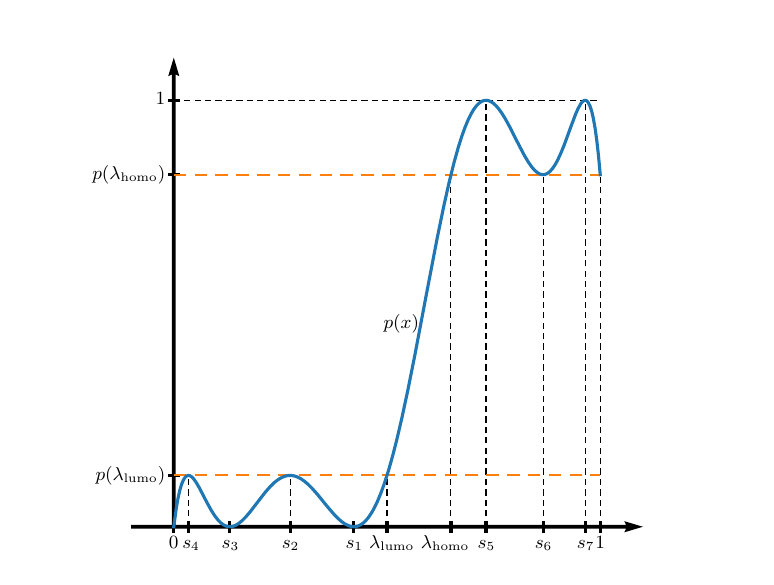}
    \caption{$L=4$, $R=3$}
    \label{fig:conditions_4_3}
  \end{subfigure}
  \hfill
  \begin{subfigure}{0.48\textwidth}
    \centering
    \includegraphics[width=\linewidth,
                     trim=40pt 10pt 60pt 15pt,clip]{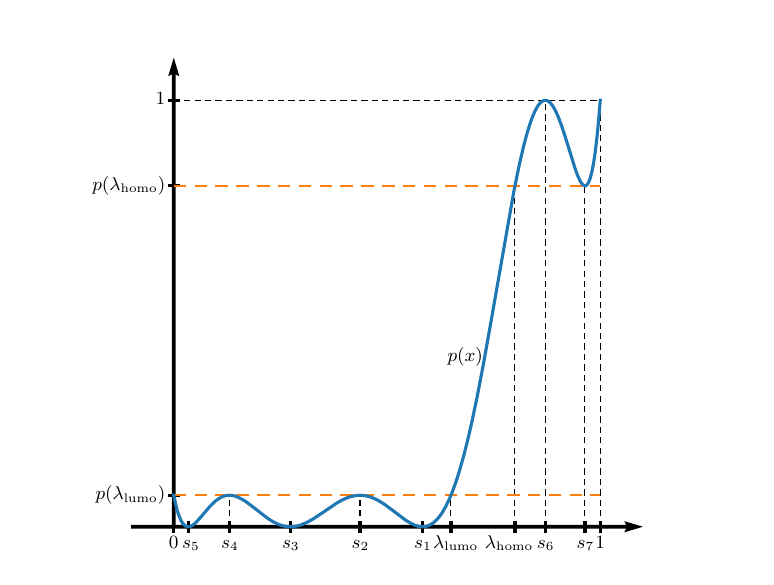}
    \caption{$L=5$, $R=2$}
    \label{fig:conditions_5_2}
  \end{subfigure}
  \caption{Illustration of the equioscillatory conditions and the resulting polynomial for two cases with different locations of $\lumo$ and $\homo$ and differing numbers of stationary points in the left and right intervals.}
  \label{fig:conditions}
\end{figure}

The equioscillatory conditions are given by 
\begin{subequations}
\label{eq:general_ansatz_conditions}
\begin{align}
  p(s_k) & =
  \begin{cases}
    0 & \textrm{if } k \textrm{ is odd}, \\
    p(\lumo) & {\rm otherwise},
  \end{cases}
  \qquad \textrm{for } k = 1,\dots,L, \label{eq:general_ansatz_conditions:left}\\
  p(s_{L+k}) & =
  \begin{cases}
    1 & \textrm{if } k \textrm{ is odd}, \\
    p(\homo) & {\rm otherwise},
  \end{cases} 
  \qquad \textrm{for } k = 1,\dots,R. \label{eq:general_ansatz_conditions:right} \\
        p(0) &=
        \begin{cases}
            0 & \textrm{if } L \textrm{ is even}, \\
            p(\lumo) & {\rm otherwise}.
        \end{cases} \label{eq:general_ansatz_conditions:zero} \\
        p(1) &=
        \begin{cases}
            1 & \textrm{if } R \textrm{ is even}, \\
            p(\homo) & {\rm otherwise}.\end{cases}\label{eq:general_ansatz_conditions:one}
    \end{align}
\end{subequations}
where
$0 \leq s_{L} \leq \dots \leq s_1 \leq \lumo < \homo \leq s_{L+1} \leq \dots \leq s_{L+R} \leq 1$.
The conditions are illustrated in Fig.~\ref{fig:conditions}.

The equioscillatory conditions give a system of nine nonlinear equations 
\begin{equation}\label{eq:nonlinsys}
    F_{(L,R;\lumo,\homo)}(s) = 0
\end{equation}
parameterized by $L$, $R$, $\lumo$, and $\homo$. For illustration, consider the case $(L,R)=(4,3)$.
The equioscillatory conditions then become the nonlinear system
\begin{align}
  F_{(4,3;\lumo,\homo)}(s) &=
  \begin{bmatrix}
    p(0) \\
    p(s_1) \\
    p(s_3) \\
    p(s_5) - 1 \\
    p(s_7) - 1 \\
    p(s_2) - p(\lumo) \\
    p(s_4) - p(\lumo) \\
    p(s_6) - p(\homo) \\
    p(1) - p(\homo)
  \end{bmatrix}
  = 0,
\end{align}
where $p(x)$ is the polynomial in \eqref{eq:polynomial_parameterization} parametrized by $s_1,\dots,s_9$ and the stationary points must satisfy the inequality constraints
\begin{equation}\label{eq:stationary-ordering}
  0 \le s_4 \le s_3 \le s_2 \le s_1 \le \lumo < \homo \le s_5 \le s_6 \le s_7 \le 1.
\end{equation}
To solve the system given in \eqref{eq:nonlinsys} we employ a
continuation-based numerical solver initialized from a set of known solutions, see Appendix~\ref{app:continuation} for details.

In the limiting case $\homo = 1$, we solve a modified problem where
\eqref{eq:general_ansatz_conditions:right} and
\eqref{eq:general_ansatz_conditions:one} are replaced by
\begin{subequations} \label{eq:modified:right}
  \begin{align}
    s_{L+r} &= 1, \quad r=1,\ldots,R,
    \\
    p(1) &= 1.
  \end{align}
\end{subequations}
Since the first $R$ derivatives vanish at $x=1$, this yields asymptotic
convergence of order $R+1$.
Similarly, when $\lumo = 0$, we
replace \eqref{eq:general_ansatz_conditions:left} and
\eqref{eq:general_ansatz_conditions:zero} by
\begin{subequations} \label{eq:modified:left}
    \begin{align}
        s_{\ell} &= 0, \quad \ell=1,\ldots,L,
        \\
    p(0) &= 0,
    \end{align}
\end{subequations}
which results in asymptotic convergence of order $L+1$ at $x=0$.
When both \eqref{eq:modified:right} and \eqref{eq:modified:left} are
imposed, the resulting polynomial is known in closed form, see
Appendix~\ref{app:continuation}, and is
monotonically increasing on $[0,\, 1]$.  

In practice, the modified
conditions are also used when $\lumo \approx 0$ and/or $\homo \approx
1$, since in this regime the oscillations provide no significant
acceleration while potentially leading to numerical instabilities when
enforcing oscillations on small intervals.
The precise criteria for switching to the modified conditions are described in the following section.

Given $\lumo$ and $\homo$, the SP8 family of polynomials is given by
\begin{align}
    \{ p_{(L,7-L;\lumo,\homo)}(x) \}_{L=0}^{7},
\end{align}
where $p_{(L,R;\lumo,\homo)}(x)$ is the unique polynomial
satisfying the equi\-oscillatory conditions for  $L$,
$R$, $\lumo$, and $\homo$.
This notation is extended to the boundary
cases $\homo=1$ and/or $\lumo=0$, where the polynomial is defined by
replacing the standard conditions with the modified conditions
\eqref{eq:modified:right} and \eqref{eq:modified:left},
respectively. In particular, $p_{(L,R;\lumo,1)}(x)$ and
$p_{(L,R;0,\homo)}(x)$ denote the corresponding one-sided boundary
cases, while $p_{(L,R;0,1)}(x)$ denotes the closed-form polynomial
obtained when both modified conditions are imposed.

\section{Recursive expansion using degree-eight polynomials}
Our recursive expansion scheme, summarized in Algorithm~\ref{alg:recursive_expansion_sp8}, builds
on the SP8 family of polynomials introduced in the previous
section. The expansion assumes the availability of upper and lower bounds for the $\lumo$ and
$\homo$ eigenvalues, satisfying the following inequalities
\begin{equation}
  0 \le \lumoouter \le \lumo \le  \lumoinner < \homoinner \le \homo \le \homoouter \le 1.
\end{equation}
Note that if the exact values are known, this can be greatly simplified by setting $\lumoinner = \lumoouter = \lumo$ and $\homoinner = \homoouter = \homo$.

\begin{figure}
  \centering  
  \begin{subfigure}{0.48\textwidth}
    \centering
    \includegraphics[width=\linewidth, trim=0pt 0pt 0pt 0pt,clip]{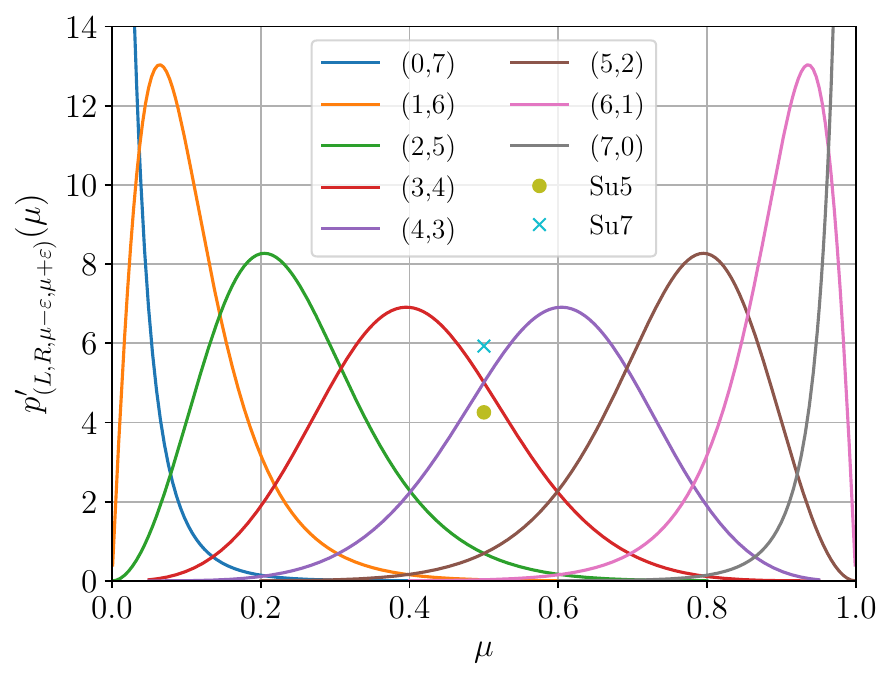}
    \caption{$p_{(L,R,\mu-\varepsilon,\mu+\varepsilon)}'(\mu)$}
    \label{fig:mu_slope}
  \end{subfigure}
  \hfill
  \begin{subfigure}{0.48\textwidth}
    \centering
    \includegraphics[width=\linewidth, trim=0pt 0pt 0pt 0pt,clip]{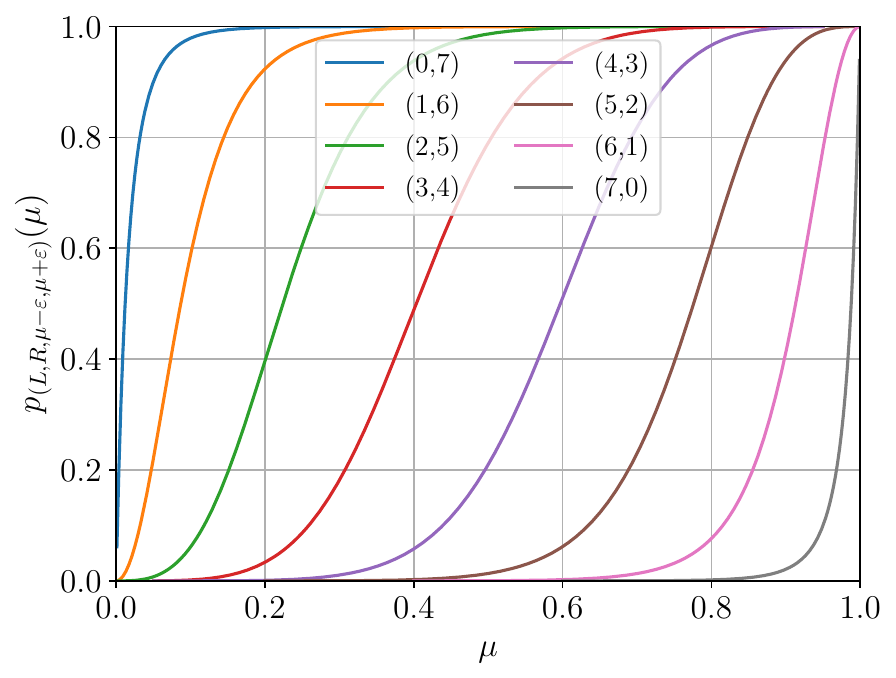}
    \caption{$p_{(L,R,\mu-\varepsilon,\mu+\varepsilon)}(\mu)$}
    \label{fig:mu_map}
  \end{subfigure}
  \caption{
For each value of $\mu$, $L$ and $R$ we define the polynomial $p_{(L,R,\mu-\varepsilon,\mu+\varepsilon)}(x)$ as the unique polynomial satisfying the oscillatory conditions in \eqref{eq:general_ansatz_conditions}. For each construction, we have used the narrow gap $\xi = 2\varepsilon,\, \varepsilon = 10^{-5}$, which is sufficiently small to represent the limiting case $\xi\rightarrow 0$. For each value of $\mu$, the corresponding polynomials and their derivatives are evaluated at $\mu$. In the figure legends, the label $(L,R)$ refers to the function $p_{(L,R,\mu-\varepsilon,\mu+\varepsilon)}( \mu )$, or its derivative.  Panel~(a) shows the gap amplification as a function of $\mu$ for each pair $(L,R)$. We have also included the gap amplification of the polynomials from the Su5 and Su7 methods for reference, using the same narrow gap. These are single points, however, since they assume $\mu=0.5$ to be fixed. Panel~(b) shows the mapping of $\mu$ for each pair $(L,R)$.
  }
  \label{fig:slope_at_mu}
\end{figure}

In each iteration of the recursive expansion we construct the SP8
family of polynomials using the outer bounds:
$\{p_{(L,7-L,\lumoouter,\homoouter)}\}_{L=0}^7$. From this family, we
select the polynomial that maximizes the gap between the transformed inner
bounds: $p(\homoinner) - p(\lumoinner)$. 
As illustrated in Fig.~\ref{fig:mu_slope}, the gap amplification varies significantly with the position of the gap.
Thus, not only does the
immediate gap amplification matter, but also the updated location of
the gap, as this determines the amplification achievable in subsequent
iterations.  Fig.~\ref{fig:mu_map} shows how $\mu$ (taken here as the center of the gap) is updated by each polynomial in the SP8 family. In particular, for gap locations that give a weaker
immediate amplification, such as gaps centered at $\mu\approx 0.3$ or
$\mu\approx 0.5$, the selected mapping tends to move $\mu$ into a
region that allows for stronger amplification in the following step.

We use a parameter $\kappa$ to determine when to deactivate acceleration and switch to the modified problems described in the previous section. Specifically, if $\lumoouter <\kappa$ we use \eqref{eq:modified:left}, and if $1-\kappa <\homoouter$ we use \eqref{eq:modified:right}. 
When the acceleration has been deactivated at both ends we start testing for convergence using a parameter-free stopping condition described in the
next section.
When the inner bounds are also both
within $\kappa$ from their target values, we start alternating between
the polynomials $p_{(4,3,0,1)}$ and $p_{(3,4,0,1)}$. Note that this gives asymptotic convergence of order $20$ in both $x=0$ and $x=1$ over two iterations.

The scheme can be run in non-accelerated mode by setting the input parameters $\lumoouter = 0$ and $\homoouter = 1$, thereby disabling acceleration from the first iteration.

\begin{algorithm}
\caption{Recursive SP8 expansion \label{alg:recursive_expansion_sp8}}
\begin{algorithmic}[1]
\Statex \textbf{Input:} $X_0$, $\lumoouter$, $\lumoinner$, $\homoinner$, $\homoouter$
\State $\kappa = 0.01$
\State $\mathbf{L}=\{0,1,\ldots,7\}$
\For{$i = 1,\ldots$}
\Statex $\quad$ // \textit{Conditional variables}
\State \makebox[2.6cm][l]{\textbf{if} $\lumoouter > \kappa$} \textbf{then} $\lumoacc = \lumoouter$ \textbf{else} $\lumoacc = 0$
\State \makebox[2.6cm][l]{\textbf{if} $\homoouter < 1- \kappa$} \textbf{then} $\homoacc = \homoouter$ \textbf{else} $\homoacc = 1$
\Statex $\quad$ // \textit{Choice of polynomial}
\State $p^{(\ell)} := p_{(\ell,7-\ell,\lumoacc,\homoacc)}$ for $\ell\in \mathbf{L}$
\State \textbf{if} $(\lumoinner < \kappa)$ \textbf{and} $(1-\kappa < \homoinner)$ \textbf{then}
\State $\quad$\textbf{if} $i$ odd \textbf{then} $\ell = 3$ \textbf{else} $\ell = 4$ // \textit{Alternate for asymptotic convergence}
\State \textbf{else}
\State $\quad$$\ell =\arg\max_{\ell \in \mathbf{L}} p^{(\ell)}(\homoinner) - p^{(\ell)}(\lumoinner)$ // \textit{Maximize gap}
\State Set: $p = p^{(\ell)}$
\Statex $\quad$ // \textit{Computation}
\State $X_i = p(X_{i-1})$
\State $(\lumoouter,\lumoinner,\homoinner,\homoouter) = (p(\lumoouter),p(\lumoinner),p(\homoinner),p(\homoouter))$
\Statex $\quad$ // \textit{Loop termination checks}
    \State \textbf{if} $\text{Tr}[X_i-X_i^2] \leq 0$ \textbf{then} \textbf{break}
\State \textbf{if} $\homoacc \neq 0$ \textbf{or} $\lumoacc \neq 1$ \textbf{then} \textbf{continue}
\State \textbf{if} $\text{Tr}[X_i-X_i^2] > C_{\text{sp8}}^{(\ell)} (\text{Tr}[X_{i-1}-X_{i-1}^2])^{q^{(\ell)}}$ \textbf{then} \textbf{break}
\EndFor
\Statex \textbf{Output:} $X_i$
\end{algorithmic}
\end{algorithm}

\subsection{Stopping criterion}\label{sec:stopping_criterion}
Following~\cite{Kruchinina_2016}, we design a stopping criterion without user-chosen tolerances that is based on monitoring discrepancies between the convergence predicted in exact arithmetic and that observed in practice.  

Given an iterative method together with an error measure, residual, or
other quantity used to assess convergence, we determine the worst-case
reduction in this quantity from one step to the next under the
assumption of exact arithmetic. For a method with asymptotic order of
convergence $q$, this reduction typically takes the form $E_{k+1} \leq
C E_k^q$. As the iteration proceeds, numerical errors, such as
floating-point roundoff, will eventually dominate, and the inequality
will no longer be satisfied.  Our parameter-free stopping criteria are
based on detecting this breakdown of the bound. Note that $C$ is
typically distinct from the asymptotic error constant of the
iteration, since the inequality is required to hold whenever
convergence is tested.

As in~\cite{Finkelstein_2021}, we assess convergence using the idempotency
measure $\textrm{Tr}[X_i-X_i^2]$. Consider an iteration of the form
$X_i = p(X_{i-1})$ with polynomial $p(x)$. Let $\{\lambda_i\}_{i=1}^n$
and $\{\gamma_i\}_{i=1}^n$ denote the eigenvalues of $X_{i-1}$ and
$X_{i}$, respectively, and assume that all eigenvalues lie in the
interval $[0,1]$.
Let $q$ be the order of the iteration associated with $p(x)$, defined as the
largest value such that
\begin{align}
  p(x)-p(x)^2 & \leq C (x-x^2)^q
\end{align}
holds for all $x \in [0,1]$ for some finite $C$. Let
\begin{align}
  C_p & = \min \{C: p(x)-p(x)^2 \leq C (x-x^2)^q \textrm{ for all } x\in[0,\, 1] \}.
\end{align}
We then have the following characterization of the convergence:
\begin{align}
  \textrm{Tr}[X_i-X_i^2] & = \sum_{i=1}^n \gamma_i - \gamma_i^2 \\
  & = \sum_{i=1}^n p(\lambda_i) - (p(\lambda_i))^2 \\
  & \leq \sum_{i=1}^n C_p (\lambda_i-\lambda_i^2)^q \\
  & \leq C_p  \Big(\sum_{i=1}^n \lambda_i-\lambda_i^2\Big)^q \\
  & = C_p \big(\textrm{Tr}[X_{i-1}-X_{i-1}^2]\big)^q.
\end{align}

The order $q$ and the constant $C_p$ depend on the polynomial that is used. 
Convergence is tested only when acceleration is not active and one of
the polynomials $p^{(\ell)} := p_{(\ell,7-\ell,0,1)}$, $\ell =
0,\ldots,7$ has been employed.
We illustrate how to determine the order $q^{(\ell)}$ and the corresponding constant $C_{\text{sp8}}^{(\ell)}$ when $\ell = 4$, i.e., for the polynomial $p^{(4)}$. Note that since the two polynomials $p^{(4)}$ and $p^{(3)}$ are related by the symmetry \eqref{eq:p_symmetry}, we have $C_{\text{sp8}}^{(4)} = C_{\text{sp8}}^{(3)}$ and $q^{(4)} = q^{(3)}$.
For $p^{(4)}$, we have
$p^{(4)}(x) = x^5(-35x^3+120x^2-140x+56)$
and
$1-p^{(4)}(x) = (1-x)^4(35x^4+20x^3+10x^2+4x+1)$.
It follows that
\begin{multline}
 p^{(4)}(x) - (p^{(4)}(x))^2  = \\
 (x-x^2)^4\underbrace{x(-35x^3+120x^2-140x+56)(35x^4+20x^3+10x^2+4x+1)}_{r(x)}
\end{multline}
which implies that $p^{(4)}(x) - (p^{(4)}(x))^2 = O((x-x^2)^4)$ as $x-x^2 \rightarrow 0$.
Since the remaining factor $r(x)$
vanishes at $x=0$ but not at $x=1$ the global order of convergence is
$q^{(4)}=4$ and the constant is given by
maximizing $r(x)$ over $[0,\, 1]$.
The factor $r(x)$ has one stationary point in $[0,\, 1]$, $x \approx
0.7578$, where it attains its maximum $C^{(4)}_{\text{sp8}} \approx 81.63$.  We set $q^{(4)}=q^{(3)}=4$ and  
$C^{(4)}_{\text{sp8}} = C^{(3)}_{\text{sp8}} = 82$ in Algorithm~\ref{alg:recursive_expansion_sp8}.
The corresponding constants for $p^{(1)}$, $p^{(2)}$, $p^{(5)}$, and $p^{(6)}$ were derived with assistance
from generative AI and thereafter verified independently. An example
prompt is included in Appendix~\ref{app:prompts}.
This yielded orders: $q^{(2)}=q^{(5)}=3$, $q^{(1)}=q^{(6)}=2$, and constants: $C^{(5)}_{\text{sp8}} = C^{(2)}_{\text{sp8}} = 56$, $C^{(6)}_{\text{sp8}} = C^{(1)}_{\text{sp8}} = 28$.
When $p^{(0)}$ or $p^{(7)}$ is applied we have no convergence at one of 0 or 1, and therefore we do not check for convergence in those cases. In Algorithm~\ref{alg:recursive_expansion_sp8} we handle this by setting $q^{(0)}=q^{(7)}=1$ and $C^{(0)}_{\text{sp8}} = C^{(7)}_{\text{sp8}} = \infty$.

\section{Numerical experiments}
We present two examples to demonstrate and evaluate the proposed recursive SP8
scheme based on polynomials of degree eight.  We first construct a
generic set of test matrices and examine how the performance depends
on two key parameters characterizing the input: the step location and the
size of the homo-lumo gap.  We then study convergence in more
detail for a Hartree-Fock calculation, in which the scheme is
used to compute the density matrix from the Fock matrix.

The new scheme is compared with the SP2 and Su5 algorithms. The SP2
and SP8 methods are considered both with and without acceleration.  Our 
implementation of the SP2 algorithm corresponds to Algorithms~4 and~5 
of~\cite{Kruchinina_2016} combined with the parameter-free
stopping criterion in~\cite{Finkelstein_2021}.  To ensure a fair
comparison, all methods are terminated using a parameter-free stopping
criterion of the same type, based on the quantity 
$\text{Tr}[X_i-X_i^2]$ as a measure of convergence, as described in
Section~\ref{sec:stopping_criterion}. 

For the Su5 method, we observed numerical instabilities for small gaps,
caused by the unstable fixed points of the component polynomial at both 0
and 1 throughout the phase in which the acceleration is active, see 
Fig.~\ref{fig:sp5_acc}. Numerical errors may cause eigenvalues to deviate
slightly outside the $[0,\, 1]$ interval, after which they drift
further away. This effect becomes more pronounced for small gaps,
where many iterations with unstable fixed points are performed. We 
therefore employ a stabilization technique that rescales the spectrum
into the $[0,\, 1]$ interval; this approach is related to a technique
used in~\cite{Amsel_2025}, adapted here to the present setting. Full details
are given in Appendix~\ref{app:su5}. No such issues were observed for 
the SP2 and SP8 schemes, with or without acceleration, where the iteration
produces polynomials that may have an unstable fixed point
at either 0 or 1, but not both.

\subsection{Generic test matrices}\label{sec:generictestmatrices}
The test matrices, used as $X_0$ in this section, are constructed for
given step location $\mu$ and gap size $\xi$. For each such pair, we
fix $n = 400$ and generate a random matrix with normally distributed
entries, from which an orthogonal matrix $Q$ is obtained via QR
factorization. The occupation number is set to $n_{\text{occ}} =
\text{round}(n(1 - \mu))$. With $\homo = \mu + \xi/2$ and $\lumo = \mu
- \xi/2$, the spectrum is constructed with $n_{\text{occ}}$
equidistant eigenvalues in $[\homo, 1]$ and $n - n_{\mathrm{occ}}$
equidistant eigenvalues in $[0, \lumo]$. The matrix $X_0$ is then
formed as $X_0 = Q\Lambda Q^T$, where $\Lambda$ is the diagonal matrix
of these eigenvalues.
For simplicity, we use $\lumoinner=\lumoouter=\lumo$ and
$\homoinner=\homoouter=\homo$ for the present test case, which corresponds to a situation where we know $\homo$ and $\lumo$ exactly. The
non-accelerated schemes are obtained by setting 
$\lumoouter = 0$ and $\homoouter=1$.

\begin{figure}[t]
  \centering
    \includegraphics[width=0.7\textwidth,trim=0pt 0pt 0pt 0pt,clip]{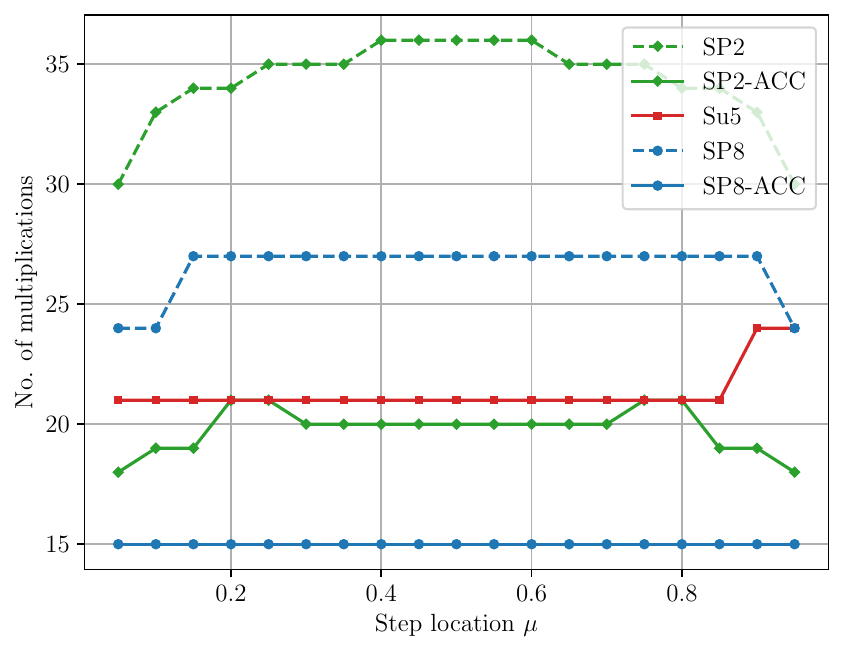}
    \caption{Number of matrix–matrix multiplications as a function of the
step location $\mu$ for the different schemes, with fixed gap size $\xi
= 0.005$. All schemes are run in double precision and terminated using
parameter-free stopping criteria as described in the main text. \label{fig:mu_vs_nmul}}  
\end{figure}

In our first test set we fix $\xi = 0.005$ and let $\mu$ vary from
$0.05$ to $0.95$. The number of matrix-matrix
multiplications needed to reach convergence as a function of the step location is shown in Fig.~\ref{fig:mu_vs_nmul}.  The accelerated SP8 scheme (SP8-ACC) outperforms all other methods across
the entire range of $\mu$, typically requiring about 25\% fewer
multiplications than the accelerated SP2 scheme (SP2-ACC), which is the fastest among the
remaining schemes, and its performance is insensitive to the step location.
The same trend is observed without acceleration,
where SP8 reduces the number of multiplications by 20-25\% compared to
SP2.

\begin{figure}[t]
    \centering
    \includegraphics[width=0.95\textwidth,trim=0pt 0pt 0pt 0pt,clip]{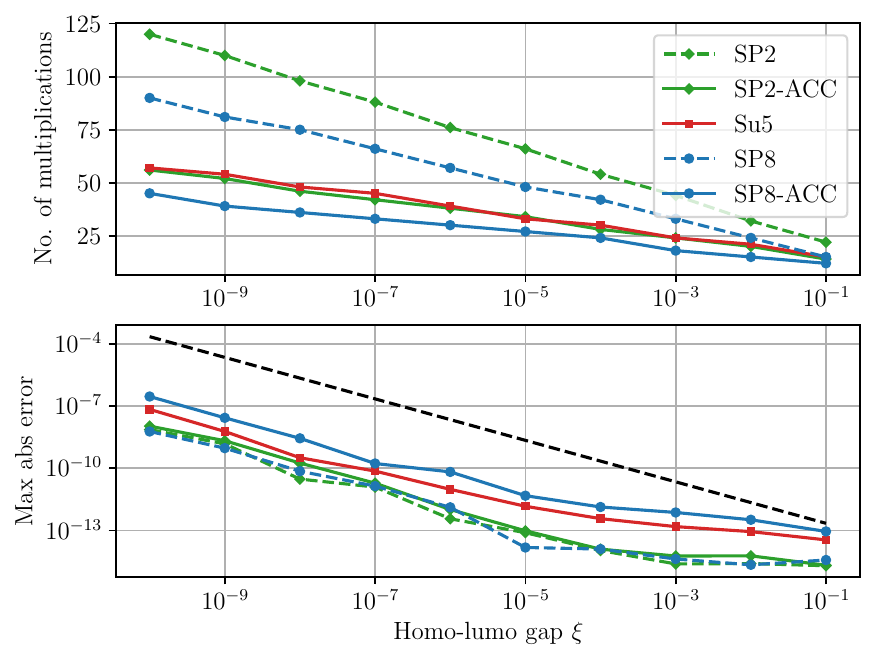}
  \caption{Number of matrix–matrix multiplications (upper panel) and maximum
error (lower panel) as functions of the homo–lumo gap size $\xi$, with
fixed step location $\mu = 0.5$. The dashed line in the lower panel
indicates $1/\xi$ scaling. All schemes are run in double precision and
terminated using parameter-free stopping criteria, as described in the main
text. \label{fig:gap_vs_nmul}}
\end{figure}

In the second test set, we instead fix $\mu=0.5$ and let $\xi$ vary
from $10^{-10}$ to $0.1$. The upper panel of Fig.~\ref{fig:gap_vs_nmul} shows the number
of multiplications as a function of the gap size, while the lower
panel shows the maximum error measured as $\max_{i,j}
|D_{i,j}-D_{i,j}^{\text{ref}}|$, where $D^{\text{ref}}$ is an accurate
reference solution. The accelerated SP8 scheme outperforms all other
methods over the entire range of gap sizes, reducing the number of
multiplications by approximately 15-25\% compared to the fastest 
competing method for each value of $\xi$.  The maximum error for
accelerated SP8 is typically one to two orders of magnitude larger
than that of the most accurate method, but exhibits the same
dependence on the condition number. The results indicate forward
stability, with the error growing roughly in proportion to the
condition number, which is $1/\xi$ as discussed in the introduction
(see~\eqref{eq:condition_number}).

\subsection{Hartree-Fock calculations}
\begin{figure}[t]
  \centering
    \includegraphics[width=0.35\textwidth,trim=0pt 70pt 0pt 50pt,clip]{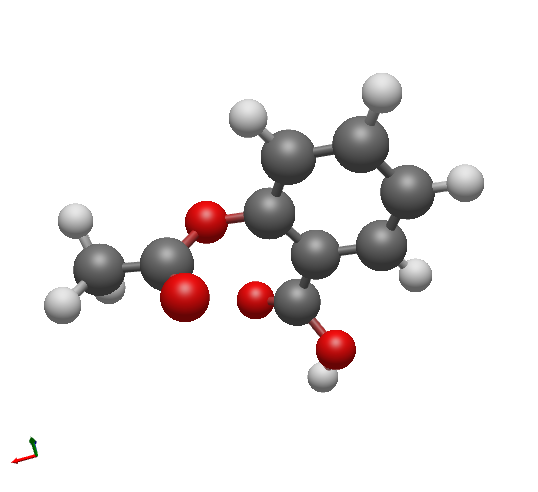}
  \caption{Ball-and-stick representation of the Aspirin molecule used in the numerical experiments.}
  \label{fig:aspirin}
\end{figure}

We consider here restricted Hartree-Fock
calculations for the Aspirin molecule, whose molecular structure is
shown in Fig.~\ref{fig:aspirin}. The molecular geometry was obtained
from PubChem (CID 2244). The Fock and overlap matrices were generated
using the Ergo quantum chemistry software, version 3.8.2~\cite{Rudberg_2018,ErgoSCF_web}, together with the standard
Gaussian basis set 6-311++G**.
The basis functions are non-orthonormal, leading to a generalized eigenvalue problem involving an overlap (Gram) matrix $S\ne I$. To transform the problem to standard form, which is assumed in Section~\ref{sec:background}, we apply a
congruence transformation $H = Z^T F Z$, where $F$ denotes the 
Fock matrix in the
non-orthonormal basis and $Z$ is an inverse factor of the overlap
matrix such that $Z^T S Z = I$. In the present calculations, the
inverse Cholesky factor of $S$ is used for this transformation.
We report results obtained from the computation of the density matrix
in the final self-consistent field iteration. 

The initial transformation in \eqref{eq:initial_transformation} is performed using extremal
eigenvalues computed with the Ergo software via the Lanczos method.
The Ergo software also computes lower and upper bounds for the
$\lumo$ and $\homo$ eigenvalues using the algorithm
of~\cite{Rubensson_2014}.  In the computations presented here, these
bounds are used in all accelerated algorithms, with the outer bounds
(i.e.\ those away from the homo-lumo gap) controlling the acceleration
in each iteration.  In the SP2 and SP8 methods, the inner
bounds are used to select the polynomial applied in each
iteration, whereas in the Su5 method only a single
polynomial is available. The spectral values are given in Table~\ref{tab:hf_spectral_values}.

\begin{table}[t]
\centering
\caption{Spectral properties and parameters for the Hartree--Fock test case.\label{tab:hf_spectral_values}}
\begin{tabular}{llllll}
  \multicolumn{6}{c}{Aspirin, HF/6-311++G**} \\
  \hline  
  $\lambda_{\min}$ & 
  $\homoouter$ & 
  $\homoinner$ & 
  $\lumoinner$ & 
  $\lumoouter$ & 
  $\lambda_{\max}$ \\
  -20.6248 &
  -0.3573 &
  -0.3482 &
  0.0355 &
  0.0565 &
  51.7152 \\
  \hline
\end{tabular}
\end{table}

\begin{figure}[t]
  \centering
    \includegraphics[width=0.95\textwidth,trim=0pt 0pt 0pt 0pt,clip]{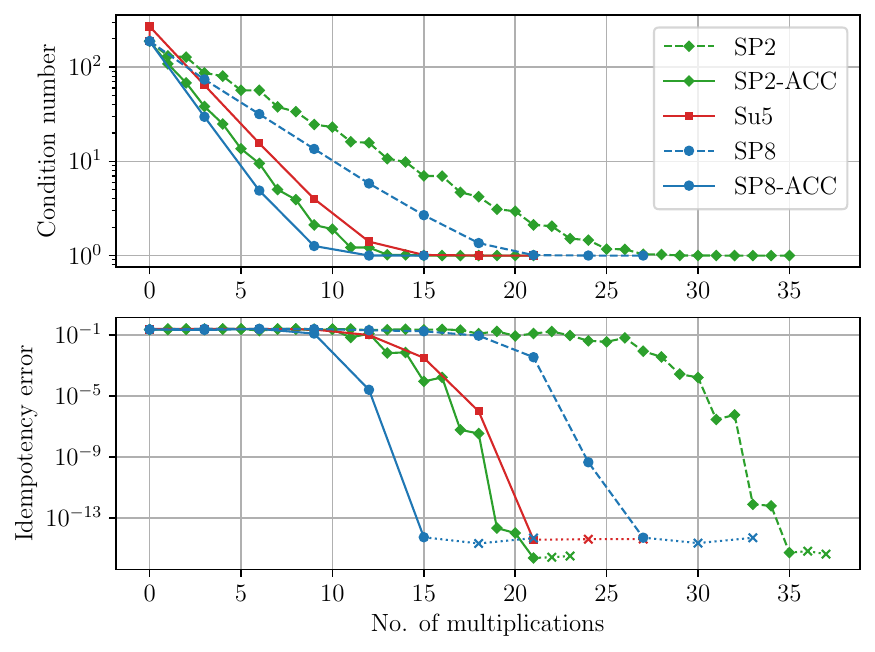}
  \caption{Condition number (upper panel, measured by $1/(\homoinner-\lumoinner)$) and idempotency error (lower panel,
measured in the max norm, $\max_{i,j}|(X-X^2)_{i,j}|$) as functions of the number of matrix-matrix
multiplications.  The test case is the computation of the density
matrix in a Hartree-Fock calculation of an aspirin molecule (see the
main text for details).  All schemes are run in double precision and
terminated using parameter-free stopping criteria.  Dotted lines with
cross markers indicate the idempotency error that would be obtained if
the iterations were continued beyond the stopping point.
  \label{fig:aspirin_convergence}}
\end{figure}

Convergence is commonly assessed using the idempotency error
$\|X_i-X_i^2\|$; however, as discussed in the introduction, this quantity only provides meaningful
information during the final purification phase of the iteration.
In the initial conditioning phase, progress is
more appropriately monitored via the condition number, given by the
inverse homo-lumo gap, see~\eqref{eq:condition_number}. The behavior of both
quantities is shown in Fig.~\ref{fig:aspirin_convergence}. During the conditioning phase, the
condition number is steadily reduced while the idempotency error
remains essentially unchanged. Once the condition number is close to
unity, the idempotency error begins to decrease rapidly, marking the
onset of the purification phase.

Consistent with Section \ref{sec:generictestmatrices}, we observe a 
significantly lower computational cost for SP8 compared to previous
methods on this real-life example. In addition, the
effectiveness of the parameter-free stopping criteria is
illustrated by plotting, using dotted lines with cross markers, the
idempotency error that would be obtained if the iterations were
continued beyond the stopping point.  The performance of the stopping
criteria is remarkable: the iterations are terminated precisely when
further iterations would no longer lead to any substantial reduction
of the idempotency error.  Accurate detection of the appropriate
termination point is particularly important for the Su5 and SP8
schemes, as each superfluous iteration incurs three additional
matrix-matrix multiplications, amounting to 20\% of the total cost for
accelerated SP8 in the present case.

\section{Discussion}
In the numerical experiments, accelerated SP8 was consistently the
fastest method to reach convergence.
As discussed in Section~\ref{sec:efficientpolynomialrepresentation}, using degree-eight component polynomials allows us to represent a larger set of polynomials compared to previous approaches, thereby enabling more efficient gap amplification.
Compared to the SP2 method, the improved performance comes at the cost of increased memory usage:
even when using the memory-efficient polynomial evaluation scheme
described in Section~\ref{sec:poly8_memory_efficient}, SP8 requires
three matrices to be stored, compared to two for SP2, highlighting a
trade-off between the number of non-scalar multiplications and memory
requirements.

In this article we measure computational cost by the number of
non-scalar multiplications, which is an accurate cost measure in the case of
large dense matrices. For sparse matrices, the situation is more
complex, as the cost of matrix-matrix multiplications depends on the
evolving sparsity pattern. The sparsity, in turn, depends on both the
polynomial sequence~\cite{Kruchinina_2016} and the scheme used to remove small matrix
elements~\cite{Rubensson_sparsity_2011}. In the sparse setting the number of non-scalar
multiplications should therefore be viewed as an approximate measure
of performance. On the other hand, the number of multiplications also
determines the critical path length, which is relevant for parallel
execution, and schemes with fewer multiplications may therefore offer
advantages in terms of parallel efficiency.  Furthermore, the number
of non-scalar multiplications is essentially independent of the
computational setup, such as hardware architecture or the choice of
linear algebra libraries, making it a reproducible performance measure.

Another trade-off arises when considering numerical accuracy.  While
all methods exhibit the expected dependence of the error on the
condition number, with the error scaling with the inverse homo-lumo
gap, we observe that, for the accelerated schemes, the error tends to
increase with the degree of the component polynomials.  This indicates
that unaccelerated methods or SP2-ACC may be more robust in regimes
with small homo-lumo gaps and/or reduced numerical precision
(e.g.\ single or lower precision arithmetic).  Alternatively, a
detailed analysis of the SP8 error may reveal how it depends on the
polynomial form and its evaluation, enabling error reduction by
adjusting the polynomial coefficients and/or the evaluation scheme.

Our recursive expansion scheme takes as input intervals containing the
$\lumo$ and $\homo$ eigenvalues. For SP2, such intervals can be computed
as a by-product of the expansion and propagated to subsequent
iterations in an outer loop, such as a self-consistent field procedure
or quantum mechanical molecular dynamics time stepping~\cite{Rubensson_2014}.  In the
present work, however, we have not addressed this aspect and instead
used either exact eigenvalues or estimates obtained from SP2 in the
numerical experiments. In future work, a similar scheme for eigenvalue
estimation could be devised for SP8.

Alternative strategies for selecting the component polynomials can be considered.  A natural approach would be to select in each iteration the degree-eight polynomial that minimizes the maximum deviation from the exact step function, i.e., to apply a greedy minimax strategy. However, such an approach turns out to be suboptimal. A consequence of choosing these polynomials is that the gap is always centered, i.e., the gap is placed at $\mu = 0.5$, but placing $\mu$ away from 0.5 allows for greater gap amplification, as illustrated in Fig.~\ref{fig:slope_at_mu}. Moreover, when $\mu=0.5$, the polynomial of degree at most eight that minimizes the maximum deviation is of degree seven, which requires four non-scalar multiplications to evaluate. In contrast, the SP8 polynomials used here remain genuinely of degree eight, since their stationary points are well separated by construction. 

We have for certain values of $\lumo$ and $\homo$ been able to construct degree-eight polynomials that amplify the gap slightly more than the SP8 polynomials over a single iteration. This means that our SP8 approach is not equivalent to a greedy gap amplification strategy, where one selects in each iteration the polynomial that maximally amplifies the gap. However, compared to other polynomials with equal or greater gap amplification that we have been able to construct, the SP8 polynomials tend to move $\mu$ to more favorable locations for subsequent iterations, making them preferable over multiple iterations. This suggests that optimal polynomial selection must account for behavior across iterations. Determining an optimal sequence of component polynomials for the recursive expansion remains an open problem, which is related to the open problems recently presented in~\cite[Section~6.3]{Amsel2026OpenQuestions}.

\appendix

\section{A continuation-based solver for SP8 polynomial parameters}\label{app:continuation}
Our solver for the nonlinear system~\eqref{eq:nonlinsys} is based on
parameter continuation and is summarized in
Algorithm~\ref{alg:parameter_continuation}.  For a given pair of
homo--lumo values $(\lumo, \homo)$, the SP8 family consists of eight
polynomials, one for each choice of $L=0,\dots,7$ stationary points in
the interval $[0, \lumo]$ and $R=7-L$ in $[\homo, 1]$.

In our algorithm, $\lumo$ and $\homo$ serve as continuation
parameters.  For each fixed pair $(L,R)$, the procedure is initialized
from a known solution $s^{\text{start}}$ corresponding to an initial
pair $(\lumo^{\text{start}}, \homo^{\text{start}})$.
In each iteration, the parameters $(\lumo^i, \homo^i)$ (with $i$ denoting the continuation step index) are moved
towards their target values $(\lumo, \homo)$ and a Newton step is
used to update the solution. For the updated parameters $(\lumo^{i+1},
\homo^{i+1})$, the Jacobian at $s^i$,
\begin{equation}
  J_{(L,R;\lumo^{i+1},\homo^{i+1})}(s^i),
\end{equation}
is computed using complex-step differentiation~\cite{Squire_1998}.

An important component of the algorithm is an adaptive step size. For
each of $\lumo^i$ and $\homo^i$, the step is restricted to a fraction
of the distance to the nearest stationary point, where the fraction is
controlled by the input parameter $\delta$. In addition, the step
selection ensures that a gap is maintained between the left and right
intervals. We used $\delta=0.5$ in all numerical experiments, which
worked well in all cases considered.

We set up the nonlinear system and corresponding initial values for
$L=4,\dots,7$, with $R=7-L$. The cases $L=0,\dots,3$ are handled using
the symmetry
\begin{align} \label{eq:p_symmetry}
  p_{(L,R,\lumo,\homo)}(x) = 1-p_{(R,L,1-\homo,1-\lumo)}(1-x)
\end{align}
which reduces the number of cases that need to be treated explicitly.

We separately treat the cases in which the conditions in
\eqref{eq:modified:right} are imposed. These give $L+2$ equations
since $R$ of the conditions are satisfied by construction. When the
conditions in \eqref{eq:modified:left} are imposed, the symmetry
relation \eqref{eq:p_symmetry} is applied as above, further reducing
the number of cases.  When both \eqref{eq:modified:right} and
\eqref{eq:modified:left} are imposed, the polynomials are known in
closed form, and no nonlinear system needs to be
solved. Table~\ref{tab:sp8_closed_form_polynomials} lists the
closed-form polynomials for $L=4,\dots,7$.  The remaining cases
$L=0,\dots,3$ can be obtained from the symmetry \eqref{eq:p_symmetry}.

\begin{table}[h]
\centering
\caption{Closed-form SP8 polynomials for $L=4,\dots,7$ (with $R=7-L$).\label{tab:sp8_closed_form_polynomials}}
\begin{tabular}{c l}
\toprule
$(L,R)$ & $p_{(L,R;0,1)}(x)$ \\
\midrule
$(4,3)$ & $-35x^8 + 120x^7 - 140x^6 + 56x^5$ \\
$(5,2)$ & $21x^8 - 48x^7 + 28x^6$ \\
$(6,1)$ & $-7x^8 + 8x^7$ \\
$(7,0)$ & $x^8$ \\
\bottomrule
\end{tabular}
\end{table}

\begin{algorithm}
\caption{Parameter continuation for the polynomial parametrization \label{alg:parameter_continuation}} 
    \begin{algorithmic}[1]
        \Statex \textbf{Input: } $L,R,\lumo,\homo, \delta$
        \State $s^0 = s^{\text{start}}(L,R)$
        \State $\lumo^0 = \lumo^{\text{start}}(L,R),\quad$ $\homo^0 = \homo^{\text{start}}(L,R)$
        \State $\xi = \homo - \lumo$
        \State $l_1 = 0,\quad$ $r_1 = 1$
        \State $i = 0$
        \While{$\lumo^i \neq \lumo$ or $\homo^i \neq \homo$}
        \State \textbf{if} $L>0$ \textbf{then} $l_1 = s_1^i$
        \State \textbf{if} $R>0$ \textbf{then} $r_1 = s_{L+1}^i$       
        \State $\Delta \lumo^{\mathrm{max}}  = \delta(\lumo^i-l_1),\quad$ $\Delta \homo^{\mathrm{max}}  = \delta(r_1 - \homo^i) $ 
        \If{$\lumo \le \lumo^i$}
        \State $\lumo^{i+1} = \max(\lumo,\lumo^i-\Delta \lumo^{\mathrm{max}})$
        \Else
        \State $\lumo^{i+1} = \min(\lumo,\lumo^i+\Delta \lumo^{\mathrm{max}},\max(\lumo^i,\homo^i-0.9\xi))$
        \EndIf
        \If{$\homo \ge \homo^i$}
        \State $\homo^{i+1} = \min(\homo, \homo^i + \Delta \homo^{\mathrm{max}})$
        \Else
        \State $\homo^{i+1} = \max(\homo, \homo^i - \Delta \homo^{\mathrm{max}},\min(\homo^i,\lumo^i+0.9\xi))$
        \EndIf
        \State $J_{(L,R;\lumo^{i+1},\homo^{i+1})}(s^i)\Delta s^i = -F_{(L,R;\lumo^{i+1},\homo^{i+1})}(s^i)$
        \State $s^{i+1} = s^i + \Delta s^i$
        \State $i = i+1$
        \EndWhile
        \State $\eta^{i-1} = \|\Delta s^{i-1}\|_{\infty}/\|s^{i}\|_{\infty}$
        \State $\eta^{i-2} = 10\eta^{i-1}$
        \While{$\eta^{i-1} > \sqrt{\varepsilon_{\!\scriptscriptstyle M}}$ or $\eta^{i-1} < \eta^{i-2}/2$}
        \State $J_{(L,R;\lumo,\homo)}(s^i)\Delta s^i = -F_{(L,R;\lumo,\homo)}(s^i)$
        \State $s^{i+1} = s^i + \Delta s^i$
        \State $\eta^i = \|\Delta s^i\|_{\infty}/\|s^{i+1}\|_{\infty}$
        \State $i = i+1$
        \EndWhile
        \Statex \textbf{Output:} $s^{i}$
    \end{algorithmic}
\end{algorithm}

\section{Prompts used to derive constants for stopping-criteria}\label{app:prompts}
We used ChatGPT~5.5 to determine the order $q$ and the constant $C_p$
needed for the parameter-free stopping criteria introduced in Section~\ref{sec:stopping_criterion}. 
We used the following query to derive the constants for iteration with the
polynomial $p_{(5,2,0,1)}$, and likewise in the other cases:
\begin{quote}
  Let
  \[
  p(x)=21x^8-48x^7+28x^6.
  \]
  Let $q$ be the order of the iteration with $p(x)$, defined as the
  largest value such that
  \[
  p(x)-p(x)^2 \leq C (x-x^2)^q
  \]
  holds for all $x \in [0,1]$ with finite
  $C$.
  Let
  \[
  C_p = \min \{C: p(x)-p(x)^2 \leq C
  (x-x^2)^q \textrm{ for all } x\in[0,\, 1] \}.
  \]
  Determine $q$ and $C_p$.  
\end{quote}

\section{The modified Su5 algorithm}\label{app:su5}
We present in Alg.~\ref{alg:recursive_expansion_su5} the modified
Su5 algorithm used in the numerical experiments,  including a
stabilization mechanism and a parameter-free stoppping criterion.

The algorithm starts with an initial linear transformation that maps
$\mu$ to $0.5$ while containing the eigenspectrum in $[0,\,1]$.  
The component polynomials $p_{(\beta)}(x)$ are polynomials of degree five
parameterized by $\beta \in [0,\, 0.5)$, determined by the conditions
$p_{(\beta)}(0)=0$, $p_{(\beta)}(x)+p_{(\beta)}(1-x)=1$, and that
$p_{(\beta)}$ equioscillates between $0$ and $p_{(\beta)}(\beta)$ at
four points on $[0,\beta]$, including the endpoints, see Fig.~\ref{fig:sp5_acc}.  For
$\beta = 0$, we define $p_{(0)}(x) = 6x^5 - 15x^4 + 10x^3$, which is
the polynomial of degree five with fixed points and two vanishing
derivatives at $0$ and $1$.

Unlike $\lumoinner$ and $\homoinner$, the outer bounds $\lumoouter$
and $\homoouter$ are not symmetric about $0.5$, yielding two possible
choices of acceleration parameter $\beta$.  On line~\ref{algline:beta_choice} of the
algorithm, the smaller value is selected to prevent immediate gap
reduction and avoid potential forward instability.  Similarly to
Alg.~\ref{alg:recursive_expansion_sp8}, we use a parameter $\kappa$ to determine when to deactivate
acceleration and instead apply $p_{(0)}(x)$.
When acceleration is active, the eigenspectrum is slightly compressed
in each iteration,
to prevent eigenvalues to drift away from $[0,\,1]$.  Without compression, such drift was
observed for small homo-lumo gaps due to the unstable fixed points of
$p_{(\beta)}$ at $0$ and $1$ for $\beta > 0$, as discussed in the main
text.  The compression is controlled by a stabilization parameter
$\sigma$.
The iteration is terminated when
\begin{equation}
  \text{Tr}[X_i-X_i^2] \leq C_{\text{su5}} (\text{Tr}[X_{i-1}-X_{i-1}^2])^3,
\end{equation}
is no longer satisfied, with $C_{\text{su5}} = 16$, derived following the approach in Section~\ref{sec:stopping_criterion} and Appendix~\ref{app:prompts}.

\begin{algorithm}
\caption{Modified Su5 expansion\label{alg:recursive_expansion_su5}}
\begin{algorithmic}[1]
  \Statex \textbf{Input:} $X_0$, $\lumoouter$, $\lumoinner$, $\homoinner$, $\homoouter$
  \State $\kappa = 0.01$
  \State $\sigma = 0.99$
\State $\mu = (\lumoinner + \homoinner) / 2$
\Statex // \textit{Choice of polynomial}
\If{$\mu \geq 0.5$}
\State Set: $p_1(x) = \frac{1}{2\mu}x$ 
\Else
\State Set: $p_1(x) = \frac{1}{2(1-\mu)}(x-1) + 1$ 
\EndIf
\State $X_1 = p_1(X_0)$
\State $(\lumoouter,\lumoinner,\homoinner,\homoouter) = (p_1(\lumoouter),p_1(\lumoinner),p_1(\homoinner),p_1(\homoouter))$
\For{$i = 2,\ldots$}
\State $\beta = \min (\lumoouter, 1-\homoouter)$ \label{algline:beta_choice}
\State \textbf{if} $\beta < \kappa$ \textbf{then} $\beta = 0$, $\sigma = 1$
\State Set: $p_i(x) = p_{(\beta)}(\sigma x +\frac{1-\sigma}{2})$
\Statex // \textit{Computation}
\State $X_i = p_i(X_{i-1})$
\State $(\lumoouter,\lumoinner,\homoinner,\homoouter) = (p_i(\lumoouter),p_i(\lumoinner),p_i(\homoinner),p_i(\homoouter))$
\Statex // \textit{Termination checks}
\State \textbf{if} $\text{Tr}[X_i-X_i^2] \leq 0$ \textbf{then} \textbf{break}
\State \textbf{if} $\beta \neq 0$ \textbf{then} \textbf{continue}
\State \textbf{if} $\text{Tr}[X_i-X_i^2] > C_{\text{su5}} (\text{Tr}[X_{i-1}-X_{i-1}^2])^3$ \textbf{then} \textbf{break}
\EndFor
\Statex \textbf{Output:} $X_i$
\end{algorithmic}
\end{algorithm}

\section*{Declarations}
As described in Appendix~\ref{app:prompts}, generative AI was used to assist in the derivation of constants appearing in the parameter-free stopping criteria. 
Generative AI was also used to assist in the writing of code for computing the coefficients of the previously established Su5 and Su7 polynomials used for comparison in this work. 
The constants and all Su5 and Su7 polynomials generated in the numerical experiments were independently verified by the authors. 

\FloatBarrier

\end{document}